\documentclass[11pt,fullpage]{article}

\usepackage{graphicx}
\usepackage{latexsym}
\usepackage{color}
\usepackage{multirow}
\usepackage{rotating}

\usepackage{amsmath}
\usepackage{graphicx}
\usepackage{amssymb}
\usepackage{amsthm}
\usepackage{amsfonts}

\renewcommand{\baselinestretch} {1.3}

\makeatletter 
\def\singlespace{\def\baselinestretch{1}\@normalsize}

\newtheorem{theorem}{Theorem}
\newtheorem{lemma}{Lemma}
\newtheorem{example}{Example}

\newtheorem{corollary}{Corollary}
\newtheorem{proposition}{Proposition}
\newtheorem{remark}{Remark}

\newcommand{\be}{\begin{equation}}
\newcommand{\ee}{\end{equation}}
\newcommand{\beqn}{\begin{eqnarray}}
\newcommand{\eeqn}{\end{eqnarray}}

\newcommand{\bt}{\begin{theorem}}
\newcommand{\et}{\end{theorem}}
\newcommand{\bl}{\begin{lemma}}
\newcommand{\el}{\end{lemma}}
\newcommand{\bp}{\begin{proposition}}
\newcommand{\ep}{\end{proposition}}
\newcommand{\bc}{\begin{corollary}}
\newcommand{\ec}{\end{corollary}}

\newcommand{\fr}[1]{(\ref{#1})} 
\newcommand{\lkr}{\left(}  
\newcommand{\rkr}{\right)} 
\newcommand{\lkv}{\left[}  
\newcommand{\rkv}{\right]} 
\newcommand{\lfi}{\left\{} 

\newcommand{\lamn}{\lambda_n}

\newcommand{\EE}{\ensuremath{{\mathbb E}}}
\newcommand{\II}{\ensuremath{{\mathbb I}}}
\newcommand{\PP}{\ensuremath{{\mathbb P}}}
\newcommand{\KK}{\ensuremath{{\mathbb K}}}

\newcommand{\Tr}{\mbox{Tr}}
\newcommand{\card}{\mbox{card\,}}

\newcommand{\bQ}{\mbox{\boldmath $Q$}}

\newcommand{\bq}{\mbox{\boldmath $q$}}

\newcommand{\bepsilon}{\mbox{\mathversion{bold}$\epsilon$}}
\newcommand{\bomega}{\mbox{\mathversion{bold}$\omega$}}

\newcommand{\supp}{\mbox{supp}\ }

\long\def\ignore#1{}

\graphicspath{%
    {converted_graphics/}
    {/}
}
\begin{document}

\title{\Large{\bf Laplace deconvolution with noisy observations}}

\author{
\large{{\sc Felix Abramovich}, {\sc Marianna Pensky}  and {\sc Yves Rozenholc}}
  \\ \\
Tel Aviv University,   University of Central Florida \\
and  Universit\'{e}  Paris Descartes  }

\date{}

\maketitle

\vspace{1cm}
\begin{abstract}
In the present paper we consider Laplace deconvolution problem for discrete
noisy data observed on an interval whose length $T_n$ may increase with the  sample size.
Although this problem arises in a variety of applications, to the best of our
knowledge, it has been given very little attention by the statistical community.
Our objective is to fill the gap and   provide   statistical analysis
of Laplace deconvolution problem with noisy discrete data.
The main contribution of the paper is   explicit construction of an
asymptotically rate-optimal (in the minimax sense) Laplace deconvolution
estimator which is adaptive to the regularity of the unknown function.
We show that the original Laplace deconvolution problem can be  reduced to
nonparametric  estimation of a regression function and its derivatives on the interval of
growing length $T_n$. Whereas the forms of the estimators   remains standard, the choices of
the parameters and the minimax convergence rates, which are expressed in terms of $T_n^2/n$
in this case,  are affected by the asymptotic growth of the length of the interval.

We derive an adaptive kernel estimator of the function of interest,
and establish its asymptotic minimaxity over a range of Sobolev classes.
We illustrate the theory by examples of construction of explicit expressions of
Laplace deconvolution estimators. A simulation study shows that, in
addition to providing
asymptotic optimality as the number of observations tends to infinity,
the proposed   estimator demonstrates good performance in finite sample examples.
 \end{abstract}


\vspace{1cm}

\bibliographystyle{plain}

\pagestyle{plain}

\section{Introduction}
\label{sec:intro}
\setcounter{equation}{0}

\subsection{Formulation and motivation}
\label{sec:form_motive}

Mathematical modeling of a variety of problems in population dynamics, mathematical physics,
theory of superfluidity and many others fields
leads to the convolution type Volterra equation of the first kind
of the form
\be \label{eq:Volterra}
q(t) = \int_0^t g(t - \tau) f(\tau) d\tau, \quad t \geq 0,
\ee
where $q(t)$ is the known  or observed  function, $g(t)$ is the  known
kernel and $f(t)$ is the unknown function to be solved for. 

Note that the LHS of equation \fr{eq:Volterra} is well defined for any $t \geq 0$
if functions $f$ and $g$ are Riemann integrable on any finite sub-interval of $[0, \infty)$.
In particular, $f$ and $g$ do not need to be absolutely or square integrable on the nonnegative
half-line. Assume the existence of their Laplace transforms
$\tilde{f}(s)$ and $\tilde{g}(s)$ for all $s \geq 0$, where 
\be \label{eq:Lapl_transform}
 \tilde{f}(s) = \int_0^\infty e^{-sx} f(x) dx,\quad  {\rm and} \quad
 \tilde{g}(s) = \int_0^\infty e^{-sx} g(x) dx,\quad s \geq 0.
\ee
In the Laplace domain,
equation  \fr{eq:Volterra} becomes $\tilde{q}(s) = \tilde{g}(s) \tilde{f}(s)$
and, therefore, the problem (\ref{eq:Volterra}) is also known as Laplace deconvolution problem.

In practice, however, one typically  has only discrete observations of the function
$q$ in \fr{eq:Volterra}  which are available only on a finite interval and, in addition,
are corrupted by noise, that leads to the following discrete noisy version of
equation \fr{eq:Volterra}
\be \label{eq:model}
y(t_i) = \int_0^{t_i} g(t_i - \tau) f(\tau) d\tau + \sigma \epsilon_i, \;\;\; i=1,...,n,
\ee
where $0 \leq t_1 \leq ... \leq t_n \leq T_n$,
$\epsilon_i$ are i.i.d. $N(0,1)$ variates,
$\sigma$ is the known constant variance  and $T_n$ may grow with $n$.

Equations of the form \fr{eq:model} appear in many practical applications.
Investigations in this paper have been motivated
by analysis of dynamic contrast enhanced imaging data and
modeling of time-resolved measurements in fluorescence spectroscopy.

\begin{example} \label{ex:DCE}
{\bf Dynamic contrast enhanced imaging data (DCE-imaging).}
\rm{DCE-imaging is widely used in cancer research
(see, e.g., Cao {\it et al.}, 2010; Goh {\it et al.}, 2005; Goh and Padhani, 2007;
Cuenod {\it et al.}, 2006; Cuenod {\it et al.}, 2011;  Miles, 2003;  Padhani and Harvey,  2005 and Bisdas
{\it et al.}, 2007). Such imaging procedures have great potential for tumor detection and
characterization, as well as for monitoring \textit{in vivo} the effects of treatments.
DCE-imaging follows the diffusion of a bolus of a contrast agent injected into a vein.
At the microscopic level, for a given unit volume voxel of interest, denote by $Y(t)$
the number of particles in the voxel   at time $t$  and by $F(t)$ the
 c.d.f. of a random lapse of time during which a particle sojourns in the voxel of interest.
Then, $F(t)$ satisfies the following equation which can be viewed as a particular case of equation \fr{eq:model}:
\be \label{eq:AIFmodel}
Y(t_i) =    \int_0^{t_i} AIF(t-\tau)(1-F(\tau)) d\tau + \sigma \epsilon_i,
\ee
where $AIF(t)$ is the Arterial Input Function which measures concentration of particles within a unit volume voxel inside 
a large artery and  can be estimated relatively easily. 
Physicians are interested in a reproducible quantification of the blood
flow inside the tissue which is characterized by $f(t)=1-F(t)$, since this quantity is independent of the number of particles of
contrast agent  injected into the vein.
}
\end{example}

\begin{example} \label{ex:fluorescene}
{\bf Time-resolved measurements in fluorescence spectroscopy.}
\rm{ Time-resolved measurements in fluorescence
spectroscopy are widely used for studies of biological
macromolecules and for cellular imaging
(see, e.g., Ameloot  and  Hendrickx, 1983;  Ameloot {\it et al.}, 1984;
Gafni, Modlin and Brand,  1975; McKinnon, Szabo and Miller, 1977;
O'Connor, Ware and Andre, 1979, and also the monograph of
Lakowicz, 2006 and references therein).
At present, in fluorescence spectroscopy, most of the time-domain measurements are carried out
using time-correlated single-photon counting. The measured intensity decay is represented
by  the number of photons $N(t_k)$  that were detected within the
time interval $(t_k, t_k + \Delta t)$, and appears as  a noisy
convolution of the  impulse response function $I(t)$ with a  known lamp function $L(t)$
$$
N(t_k) = \int_0^{t_k} L(t_k -\tau) I(\tau) d\tau + \sigma \epsilon_k.
$$
The objective is to determine the impulse response function
$I(x)$  that best matches the experimental data.
}
\end{example}

\subsection{Difficulty of the problem}
\label{sec:difficulty}

The mathematical theory of (noiseless) convolution type Volterra equations
is well developed (see, e.g., Gripenberg, Londen and Staffans, 1990)
and the exact solution of \fr{eq:Volterra} can be obtained through
Laplace transform.  However, direct application of Laplace transform for
discrete measurements faces serious conceptual and numerical problems.
The inverse Laplace transform is usually  found by application of tables of inverse Laplace
transforms, partial fraction decomposition or
series expansion (see, e.g., Polyanin and Manzhirov, 1998), neither of which is
applicable in the case of the discrete noisy version of Laplace deconvolution.

Formally, by extending $g(t)$ and $f(t)$ to the negative values of $t$ by setting
$f(t) = g(t)= 0$ for $t<0$, equation
\fr{eq:Volterra} can be viewed as a particular case of the
Fredholm convolution equation
\be \label{eq:Fredholm}
h(t) = \int_{-\infty}^\infty  g(t  - \tau) f(\tau) d\tau,
\ee
whose discrete stochastic version
\be \label{eq:fourier_deconv}
y(t_i) = \int_{-\infty}^\infty  g(t_i - \tau) f(\tau) d\tau + \sigma \epsilon_i, \;\;\;\; i=1,...,n,
\ee
known also as Fourier deconvolution problem,
has been extensively studied in the last thirty years
(see, for example, Carroll and Hall, 1988; Comte, Rozenholc and Taupin, 2006,
2007; Delaigle, Hall and Meister, 2008;
Diggle and Hall, 1993; Fan, 1991; Fan and Koo, 2002; Johnstone {\it et al.},
2004; Pensky and Vidakovic, 1999;  Stefanski and Carrol, 1990 among others; see also  monograph
by Meister,2009 and references therein).

Unfortunately, the existing approaches to Fourier deconvolution cannot be easily
extended to solution of noisy discrete version of
Laplace convolution equation \fr{eq:model}.
The body of work cited above addresses one of three
situations:  the case when functions $f$ and $g$ are periodic with
period $T$, density deconvolution, and the case of random design, where
$t_i$ in \fr{eq:fourier_deconv} are random variables generated by
some density function.

In the first setup, convolution
\fr{eq:Fredholm} becomes circular convolution and measurements in
equation \fr{eq:fourier_deconv} are taken on an interval of fixed
length $T$, so that the problem can be solved by application of
discrete Fourier transform.
However, since the functions $f$ and $g$ are not periodic on $[0, T_n]$, the
integral in the RHS of
equation \fr{eq:model} is not a circular convolution and the discrete Fourier
transform
cannot be directly applied.
Furthermore, the length of the interval $T_n$ may grow with $n$
that affects the convergence rates. For relatively small $T_n$
(e.g., $T_n \sim \log n$), approximation
of Fourier transform by its discrete version will be very poor
which results in low convergence rates of the estimator of $f$.

Density deconvolution problem and
nonparametric regression estimation with random measurements $t_i$
typically assume that, as $n \to \infty$, the measurements $t_i$  in
\fr{eq:fourier_deconv} adequately represent the domain of $h(t)$ in
\fr{eq:Fredholm}.
In these setups observations are absent on a particular
part of the domain only if the density which generates those observations is very low.
This, however, is not at all true for equation \fr{eq:model}  where
lack of observations for $t>T_n$ is due entirely to experimental design and has
no relation to the values of the estimated function.

To the best of our knowledge, nobody tackled the
problem of Fourier deconvolution \fr{eq:fourier_deconv} when
observations $t_i$ are non-random fixed quantities on an interval of
 length $T_n$ which grows with the number of observations.
In addition, we should also mention the important {\em causality} property
of the Laplace deconvolution not shared by its Fourier counterpart, where the
values of $q(t)$ for $0\leq t \leq T_n$
depend on values of $f(t)$ for $0\leq t \leq T_n$
only and vice versa.
Finally, we show that under mild conditions, the solution of the equation \fr{eq:model}
can be represented {\em explicitly} via derivatives of the RHS $q(t)$ that implies computational advantages of the proposed approach.

\subsection{Existing results}
\label{sec:bibliography}

Only few applied mathematicians took an effort to tackle the problem with
discrete measurements in the LHS of \fr{eq:Volterra}.
Ameloot and Hendrickx  (1983) applied Laplace deconvolution for the
analysis of  fluorescence curves and used a parametric presentation of the solution $f$ as a sum of exponential functions
with parameters evaluated  by minimizing discrepancy with the RHS.
In a somewhat similar manner,  Maleknejad {\it et al.} (2007) proposed to expand
the unknown solution over a wavelet basis and
find the coefficients via the least squares algorithm.  Lien {\it et al.} (2008), following Weeks  (1966),
studied numerical inversion of the Laplace transform using Laguerre functions.
Finally, Lamm (1996) and Cinzori and Lamm (2000) used discretization of
the equation \fr{eq:Volterra} and  applied various versions of the Tikhonov regularization  technique.
However, in all of the above papers, the noise in the measurements  was either ignored
or treated as deterministic. The presence of random noise in (\ref{eq:model})
makes the problem even more challenging.

Unlike Fourier deconvolution that has been intensively studied in statistical
literature (see references above), Laplace deconvolution
received virtually   no attention within statistical framework.
To the best of our knowledge, the only paper which tackles the problem is
Dey, Martin and Ruymgaart (1998) which considers a noisy version of Laplace deconvolution
with a very specific kernel of the form $g(t)=be^{-at}$. The authors use the fact that, in this case,
the solution of the equation \fr{eq:Volterra}  satisfies a particular linear differential equation and,
hence, can be recovered using   $q(t)$ and its derivative $q' (t)$.
For this particular kind of kernel, the authors   derived   convergence
rates for the  quadratic risk of the proposed estimators, as $n$ increases,
under the assumption that  the $m$-th derivative of
$f$ is   continuous on $(0,\infty)$. However, they assume
that data is available on the whole nonnegative half-line (i.e. $T_n = \infty$)
and that $m$ is known (i.e., the estimator is not adaptive).

\subsection{Objectives and organization of the paper}
\label{sec:objectives}

For   the  reasons listed above, estimation of $f$ from discrete noisy
observations $y$ in \fr{eq:model} requires development of a novel approach.
The objective of the present paper is to fill the gap and to develop general
statistical methodology for Laplace deconvolution problem which allows to circumvent lack
of observations for $t> T_n$ and leads to effective representation of $f$ on the interval $(0,T_n)$,
no matter what value $T_n$ takes.
We establish minimax convergence rates for Laplace deconvolution setup over
Sobolev classes and
derive the adaptive estimator of $f$ which is rate-optimal
over entire range of Sobolev classes. The proposed estimator
is based on estimating $q$ and its derivatives from noisy data $y$ in 
\fr{eq:model}, where $q(t)=(f \ast g)(t)=\int_0^t g(t-\tau)f(\tau)d\tau$ is the
convolution of $f$ and $g$.
Thus, one can use the numerous existing techniques for nonparametric estimation of a function and its derivatives.
In particular, we employ kernel estimators with the global bandwidth adaptively
selected by Lepski procedure.

An attractive feature of the estimation technique proposed in this paper is that estimator
of $f$ is expressed \emph{explicitly} via $q$ 
and its derivatives.
Another interesting aspect of the  considered model (\ref{eq:model}) is
that the data is observed on the interval of asymptotically increasing length, where
$T_n \rightarrow \infty$   as $n \rightarrow \infty$.
This is indeed a reasonable assumption since, as $n$ is growing, demands on the improvements
of the estimation precision require to decrease the bias by sampling $q(t)$ for larger and larger  values of $t$.
Dependence of $T$ on $n$ may not significantly affect  estimation procedures but evidently leads to different convergence rates
that are formulated in terms  of $T_n^2/n$.

The rest of the paper is organized as follows.
Section~\ref{sec:main} delivers main results of the paper. In particular,
Section~\ref{sec:assumptions} introduces notations and assumptions used
throughout the paper. In Section~\ref{sec:lower_bounds} we derive the lower bounds on the minimax risk
of estimating $f$ in \fr{eq:model}. Section~\ref{subsec:explicit}
reviews some mathematical results for noiseless Laplace deconvolution relevant
for constructing the proposed estimator.
Section~\ref{subsec:noisy_solution} is dedicated to explicit derivation of
Laplace deconvolution estimator in model \fr{eq:model}, while Section~\ref{sec:adaptive_minimax}
establishes its asymptotic adaptive minimaxity over entire range of Sobolev classes.
Section~\ref{sec:examples} contains examples of explicit estimators of Laplace deconvolution for various
types of  kernels $g$. The results of a simulation study are presented in
Section~\ref{sec:simul}. Section~\ref{sec:discussion} concludes the paper with discussion.
All the proofs are given in Appendix.


\section{Main results}
\label{sec:main}
\setcounter{equation}{0}


\subsection{Notations and assumptions}
\label{sec:assumptions}
\setcounter{equation}{0}

In this section we introduce   notations and assumptions used throughout the paper.

The $L_k(\mathbb{R}^+)$-norm of the function $h$ is denoted  by $\|h\|_k$ and
$\|h\|_\infty$ is the supremum norm of $h$.
If $k=2$ and there is no ambiguity, we shall omit the subscript
in the notation of the norm, i.e.   $\|h\| = \|h\|_2$.
We use the standard notation $W^{r,p}(\mathbb{R}^+)$ for a Sobolev space of functions on $[0,\infty)$ that
have $r$ weak derivatives with finite $L_p$-norms and omit $p$ in this notation  if $p=2$, that is,   $W^r(\mathbb{R}^+)=W^{r,2}(\mathbb{R}^+)$.
In addition, we shall omit $\mathbb{R}^+$ in the notations of the norms and functional spaces and,
unless the opposite is stated, assume
that all functions are defined on the nonnegative part of the real line.
\\

Let $r \geq 1$  be such that
\be \label{k_cond}
g^{(j)} (0) = \lfi
\begin{array}{ll}
0, & \mbox{if}\ \ j=0, ..., r-2,\\
B_r \ne 0, &  \mbox{if}\ \ j=r-1,
\end{array} \right.
\ee
with obvious modification $g(0)=B_1 \ne 0$ for $r=1$.

Assume now the following conditions on the unknown $f$ and the  known kernel
$g$ in (\ref{eq:Volterra}):

\begin{itemize}
\item
[(A1)]  $g \in W^{r,1} \cap W^{\nu}$, $\nu \geq r$.

\item
[(A2)] Let $\Omega$ be a collection of distinct zeros $s_\omega$ of the
Laplace transform $\tilde{g}$ of $g$.
Then all zeros of $\tilde{g}$ have negative real parts, i.e.,
$$
s^* = \max_{s_\omega  \in \Omega}  Re (s_\omega)  < 0.
$$

\item
 [(A3)] $f \in W^{m}$ where $m \leq \nu +1 -r$.
\end{itemize}

\ignore{
\begin{remark} \label{rm:assump}
\rm{
 In many cases, $\nu = \infty$ and there are no restrictions on the smoothness
of $f$. On the other hand, if $\nu$ is finite, there is no advantage in taking $m$ larger than $\nu +1-r$ since
it will not translate into the better convergence rates. It is very easy to explain this
phenomenon at an intuitive level. If $f(t)$ and $g(t)$ are, respectively, s $m$ and $\nu$ times
continuously differentiable for $t>0$, then $q$ is $L$ times continuously differentiable where
$L = \min(r+m+1, \nu + m)$, so that its $r$-th derivative is $L-r$ times continuously differentiable.
Since estimator of $f$ is based on estimating $q^{(r)}$, one needs  $q^{(r)} \in W^{m}$ to attain adaptivity
which leads to the condition $m \leq \nu +1 -r$.  If one works with Sobolev spaces,
this imposes natural restriction on the range of adaptivity. However, since $g$ can have only
jump discontinuities in its derivatives, consideration of, e.g., Besov spaces can extend adaptivity to a wider range of spaces.
}
\end{remark}
}

Finally, we impose the following assumption on $T_n$ and design points $t_i$, $i=1,..., n$:

\begin{itemize}
\item[(A4)]  Let $T_n$ be such that $T_n \rightarrow \infty$ but $n^{-1}\, T_n^2 \rightarrow 0$ as $n \rightarrow \infty$
and there exist $1 \leq \mu < \infty$ such that  $\max_i\, |t_i-t_{i-1}| \leq \mu  n^{-1}T_n$.
\end{itemize}

In what follows, we use the symbol $C$ for a generic positive constant, independent
of the sample size $n$, which may take different values at different places.


\subsection{Lower bounds for the minimax risk}
\label{sec:lower_bounds}

In order to establish a benchmark for an estimator of an unknown function $f$ from
its noisy Laplace convolution \fr{eq:model} we derive the asymptotic
minimax lower bounds for the $L_2([0,T_n])$-risk over a Sobolev ball
$W^{m}(A)$ of radius $A$.
It turns out that, unlike in the density deconvolution problem or Fourier deconvolution setup,
the rates of convergence depend on the length of the interval $T_n$ and are expressed in terms of the ratio $T_n^2/n$:

\begin{theorem} \label{th:lower_bounds}
Let condition \fr{k_cond} and Assumptions   (A1)--(A4) hold.  Then, there exists a constant $C>0$ such that
\be \label{eq:lower_bound}
\inf_{\hat{f}_n} \sup_{f \in W^m(A)} E||\hat{f}_n-f||^2_{L_2([0,T_n])} \geq
C \lkr \frac{T_n^2}{n} \rkr^\frac{2m}{2(m+r)+1},
\ee
where the infimum is taken over all possible estimators $\hat{f}_n$ of $f$, and, therefore,
%
$$
\inf_{\hat{f}_n} \sup_{f \in W^m(A)} E||\hat{f}_n-f||^2_{L_2([0,\infty))} \geq C \lkr  \frac{T_n^2}{n} \rkr^\frac{2m}{2(m+r)+1}.
$$
\end{theorem}


\subsection{Solution of noiseless Volterra equation}
\label{subsec:explicit}

As we have already mentioned, unlike Fourier deconvolution, an estimator $\hat{f}_n$ of the unknown $f$ in
\fr{eq:model} can be obtained explicitly in the closed form.
To understand the motivation for the proposed $\hat{f}_n$ we find first the
exact solution of the noiseless Volterra equation
(\ref{eq:Volterra}).

Taking derivatives of both sides of
\fr{eq:Volterra} under \fr{k_cond} and Assumptions (A1), (A3), one obtains
\begin{eqnarray}
q^{(j)} (t) &=& \int_0^t g^{(j)}(t-\tau) f(\tau) d\tau,\ \ j=1,..., r-1;
\nonumber\\
\cdot & \cdot & \cdot \nonumber\\
q^{(r)} (t) &=& B_r f(t) + \int_0^t g^{(r)}(t-\tau) f(\tau) d\tau,
\label{eq:qr}
\end{eqnarray}
which is the Volterra equation of the second kind.
Taking higher-order derivatives, (\ref{eq:qr}) yields
\begin{eqnarray*}
q^{(r+1)}(t)&=&B_r f'(t)+g^{(r)}(t)f(0)+\int_0^t g^{(r)}(t-\tau)f'(\tau)d\tau, \\
\cdot & \cdot & \cdot \\
q^{(r+m)}(t)&=&B_r f^{(m)}(t)+\sum_{j=0}^{m-1}g^{(r+j)}(t)f^{(j)}(0)+
\int_0^t g^{(r)}(t-\tau)f^{(m)}(\tau)d\tau.
\end{eqnarray*}
Then, under Assumptions (A1) and (A3), one has $q^{(r+m)} \in L_2$ and,
hence, $q \in W^{r+m}$.

In addition, due to Assumptions (A1) and (A3), condition (\ref{eq:qr}) implies that
$q^{(r)} \in L_1$ and, therefore, one can use
the following known facts from the theory of Volterra equations of the second
kind:
\begin{enumerate}
\item there exists a unique solution $\phi$ of the equation
\begin{equation} \label{eq:new_eq}
g^{(r)}(t) = B_r \phi(t) + \int_0^t g^{(r)}(t-\tau) \phi(\tau) d\tau
\end{equation}
called    a {\em resolvent} of $g^{(r)}$
(see Theorem 3.1 of Gripenberg, Londen and Staffans, 1990);
\item
there exists a unique solution $f$ of \fr{eq:qr} and, therefore, of the
original equation \fr{eq:Volterra}, which can be written as
\be \label{eq:solution_f}
f(t) = B_r^{-1} q^{(r)} (t) - B_r^{-1} \int_0^t q^{(r)}(t-\tau) \phi(\tau) d\tau
\ee
(see Theorem 3.5 of Gripenberg, Londen and Staffans, 1990).
\end{enumerate}

\begin{remark}  \label{rem:assumA2}
\rm{Assumption (A2) ensures that the solution $f(t)$ of the noiseless
Laplace convolution equation (\ref{eq:Volterra}) is numerically stable.
By Half-Line Paley-Wiener theorem (see Theorem 2.4.1 of Gripenberg, Londen and Staffans, 1990),
the resolvent $\phi(\tau)$ of $g$ in \fr{eq:new_eq} is absolutely integrable if and only if Assumption (A2) is satisfied. If
$\tilde{g}$ has roots with positive real parts, then, by Corollary 2.4.2 from 
the same book, $\phi(\tau)$ is growing at an
exponential rate, so that
$e^{-s \tau} \phi(\tau)$ is absolutely integrable for any $s > s^*$,
where $s^*$ is defined in assumption (A2).}
\end{remark}

It follows from the above  that, in order to solve the noiseless Volterra equation
(\ref{eq:Volterra}), one only needs to determine a resolvent
$\phi$ in (\ref{eq:new_eq}) defined entirely by the $r$-th derivative $g^{(r)}$
of the ({\em known}) kernel $g$.
Taking Laplace transform of both sides of (\ref{eq:new_eq}) yields
$$
\widetilde{g^{(r)}}(s)=B_r\tilde{\phi}(s)+\widetilde{g^{(r)}}(s)\tilde{\phi}(s)
$$
where, due to \fr{k_cond}, one has $\widetilde{g^{(r)}}(s)=s^r \tilde{g}(s)-B_r$.
Therefore, $\phi(t)$ can be obtained as an inverse Laplace transform
of $\tilde{\phi}$, where
\be \label{eq:Lapl_phi}
\tilde{\phi}(s)=\frac{s^r\tilde{g}(s)-B_r}{s^r\tilde{g}(s)}.
\ee

Behavior of the resolvent function $\phi$ is thus determined by the properties
$\tilde{g}$.
It turns out (see, e.g., Gripenberg, Londen and Steffans 1990, Chapter~7) that,
under Assumption  (A2) and \fr{k_cond}, $\tilde{g}$ is analytic  and, hence,
all its zeros are well separated. Moreover, $\phi$ can be presented as the sum of a polynomial
of degree $(r-1)$ and an absolutely integrable function.
In a variety of practical applications, the kernel $g$ is represented by a
combination of some elementary functions and, hence, $\tilde{g}$ is not an oscillating function.
 Hence,   the number of zeros of $\tilde{g}$ is finite and, since $\tilde{g}$ is an analytic function,
these zeros are of finite orders. In this case, solution $f$ can be written explicitly
as it  follows from the following theorem:

\begin{theorem}  \label{th:semiexplicit}
Let condition \fr{k_cond} and Assumptions (A1)--(A3) hold.
Then, the resolvent $\phi$ in \fr{eq:Lapl_phi}  is of the form
\be \label{eq:phi_form1}
\phi(t) =   \sum_{j=0}^{r-1} \frac{a_{0,j}}{j!} t^j   + \phi_1 (t),
\ee
where $\phi_1 \in L_1$. Hence, by \fr{eq:solution_f},
$f$ in  \fr{eq:Volterra} can be recovered as
\be \label{eq:f_semiexplicit}
f(t) = B_r^{-1}\left(q^{(r)} (t) - \sum_{j=0}^{r-1} a_{0, r-1-j} q^{(j)} (t)  -
\int_0^t q^{(r)}(t-\tau) \phi_1 (\tau) d\tau \right).
\ee
If, in addition, $\tilde{g}$ has a finite number $M$ of  distinct zeros of
orders $\alpha_l$, respectively, $l=1, ..., M$, then  $f$ is of the form
\be \label{eq:f_explicit}
f(t) =  B_r^{-1} \left( q^{(r)} (t) -  \sum_{j=0}^{r-1} b_j q^{(r-1-j)} (t)
- \int_0^t q(t-x) \phi_1^{(r)}(x) d x
\right),
\ee
where $s_0=0$, $\alpha_0 =r$ and
\begin{eqnarray}
\phi_1(x)    & = &   \sum_{l=1}^M \sum_{j=0}^{\alpha_l -1} \frac{a_{l,j} x^j e^{s_l x}}{j!},\label{eq:psi}\\
a_{l,j}  & = &  \frac{1}{(\alpha_l -1-j)!}\  \frac{d^{\alpha_l -j-1}}{ds^{\alpha_l -j-1}} \lkv (s - s_l)^{\alpha_l}
\tilde{\phi}(s) \rkv \Bigg|_{s=s_l}\;,  \label{eq:coefs}\\
b_j   & = &  a_{0,j} + \sum_{l=1}^M \sum_{i=0}^{\min(j, \alpha_l-1)} {j \choose i} a_{l,i} s_l^{j-i}\;. \label{eq:bj}
\end{eqnarray}
\end{theorem}

\begin{remark} \label{rem:zeros}
\rm{
Note that in Theorem \ref{th:semiexplicit},
Assumption (A1)  and condition \fr{k_cond}   are essential for explicit construction of estimators.
However, calculations  in \fr{eq:phi_form1}--\fr{eq:bj} can be carried out without Assumption (A2) being valid.
Assumption (A2) is only needed to ensure that $\phi_1 \in L_1$. In particular, if the number of zeros is finite,
then  $Re(s_l)<0$, $l=1,...,M$, implies that $\phi_1$ in \fr{eq:psi} is a sum of products of polynomials and
exponentials with powers having negative real parts and, hence,
$\phi_1 \in L_1 \cap L_2$. If some of zeros have positive real parts,
expansions \fr{eq:coefs} and \fr{eq:bj} in Theorem \ref{th:semiexplicit} will still be valid
but $\phi_1^{(r)}$ will contain exponential terms with positive powers that will
grow and magnify the errors of estimating $q$ as $t$ tends to infinity.
}
\end{remark}


\subsection{Adaptive estimation of Laplace deconvolution}
\label{subsec:noisy_solution}
Theorem \ref{th:semiexplicit}  leads to  an
estimator $\hat{f}_n$ in \fr{eq:model}  of the  semi-explicit form
\be \label{eq:estimator}
\hat{f}_n(t)=B_r^{-1} \lkr \widehat{q^{(r)}} (t) -
\sum_{j=0}^{r-1} a_{0,r-1-j} \widehat{q^{(j)}} (t)  -
\int_0^t \widehat{q^{(r)}}  (t-\tau) \phi_1 (\tau) d\tau \rkr,
\ee
where $\widehat{q^{(j)}}(t)$ are some estimators of $q^{(j)} (t)$,
$j=0,\ldots,r$, and
the function $\phi_1$ is expressed in terms  of the inverse Laplace transform
of the completely known function  $\tilde{\phi}$ defined in \fr{eq:Lapl_phi}.
Under the additional (usually satisfied) condition that $\tilde{g}$ has a finite number of zeros, the second statement of
Theorem \ref{th:semiexplicit} leads to an explicit expression for the estimator with $\phi_1$ defined by \fr{eq:psi}:
\be \label{eq:f_explicit_est}
\hat{f}_n(t) =  B_r^{-1} \left( \widehat{q^{(r)}} (t) -  \sum_{j=0}^{r-1} b_j \widehat{q^{(r-1-j)}} (t)
- \int_0^t \widehat{q}(t-x) \phi_1^{(r)}(x) d x.
\right)
\ee

Note that, unlike \fr{eq:estimator}, the integral term in
\fr{eq:f_explicit_est} involves $q$ rather than $q^{(r)}$ and, hence,
the boundary effects of estimating derivatives do not propagate to interior
points of the interval $[0,T_n]$.

Laplace deconvolution can be therefore   reduced to nonparametric
estimation of $q=f \ast g \in W^{r+m}$ (see Section \ref{subsec:explicit})
and its derivatives of orders up to $r$ from the discrete
noisy data in the model
$$
y(t_i)=q(t_i)+\sigma \epsilon_i,\quad i=1,\ldots,n,
$$
where $0 \leq t_1 \leq \ldots \leq t_n \leq T_n$, $\epsilon_i$ are i.i.d.
$N(0,1)$ variates and $\sigma>0$ is known.
This is a well-studied problem, and estimation can be carried out  by
a number of various approaches, e.g.,
kernel estimation, splines, local polynomials, wavelets, etc.

It is important to note however that for the problem at hand, the data is sampled on
an interval of asymptotically increasing length that calls for
necessary modifications of traditional estimators and affects their global
convergence rates on the interval $[0,T_n]$ which are   expressed in terms of $n^{-1} T_n^2$.

For illustration, we consider kernel estimation with the global bandwidth
selected adaptively by Lepski technique.
To estimate the $j$-th derivative of $q(t)$, $j=0,\ldots, r$, $t \in [0, T_n]$,
choose a kernel function $K_j$ (not to be confused with the convolution kernel $g$) of order $(L,j)$ with $L>r$
satisfying the following conditions:
\begin{enumerate}
\item[(K1)] $\supp(K_j) = [-1,1]$, $K_j$ is twice continuously differentiable
and $\int K_j^2(t)dt < \infty$.
\item[(K2)] $\int t^l K_j(t)dt=
     \left\{
        \begin{array}{ll}
         0, & l=0,...,j-1,j+1,...,L-1, \\
         (-1)^j j!, & l=j.
        \end{array}
     \right. $
\end{enumerate}
Construction of such kernels is described in, e.g., Gasser, M\"uller and Mammitzsch (1985).

Define a well-known Priestley-Chao type kernel estimator of $q^{(j)}$
with a (global) bandwidth $\lambda_j$:
\be \label{eq:pc}
\widehat{q_{\lambda}^{(j)}}(t)= \frac{1}{\lambda_j^{j+1}}\sum_{i=1}^n
K_j\left(\frac{t-t_i}{\lambda_j}\right)(t_i-t_{i-1})y_i.
\ee
Certain routine boundary corrections are required for $t$ close to the boundaries
(see Gasser and M\"uller, 1984 for details).

We utilize a general methodology
developed by Lepski (e.g., Lepski, 1991) for data-driven selection
of a bandwidth $\lambda_j$ in \fr{eq:pc}.
In particular, we apply the global bandwidth version of Lepski, Mammen and
Spokoiny's (1997) procedure and modify it also for estimating derivatives.

The resulting procedure for choosing $\lambda_j$ in \fr{eq:pc} can be described as follows.
For each $j$, $0 \leq j \leq r$, and the corresponding kernel $K_j$ of order
$(L,j),\;L>r$, consider
the geometric grid of bandwidths $\Lambda_j$, where
\be \label{LambdaDef}
\Lambda_j =\{ \lambda_l=  a^{-l},\;l=0,1,..., J_n;\ J_n =(2j+1)^{-1}\, \log_a (n\, \sigma^{-2} \, T_n^{-2})  \},
\ee
and $a>1$ is an arbitrary constant.
Smaller values of $a$ allow a finer choice of the optimal bandwidth but
increase computational complexity. 
Note that cardinality of $\Lambda_j$ does not exceed $\log_a n$ since
$\card (\Lambda_j)  = 1 + J_n  \leq \log_a n$.
Define
\be \label{eq:adapband}
\hat{\lambda}_{j,n} =\max\left\{\lambda \in \Lambda_j: \|\widehat{q^{(j)}_\lambda}-\widehat{q^{(j)}_h}\|^2_{[0,T_n]}
\leq \frac{ 4\,  C_j^2\, \sigma^2 T_n^2}{ n\, h^{2j+1}}\;
{\rm for\; all\;} h \in \Lambda_j,  \left(\frac{\sigma^2  T_n^2 }{n}\right)^{\frac{1}{2j+1}} \leq  h < \lambda  \right\},
\ee
where constants $C_j$ are such that
\be \label{c0}
C_j^2 >   \mu^2  \|K_j\|^2
\ee
and $\mu$ is defined in Assumption (A4).

We then estimate $q^{(j)}$ by
\begin{equation} \label{eq:deriv_est}
\widehat{q_{\hat{\lambda}_{j,n}}^{(j)}}(t)= \frac{1}{\hat{\lambda}_{j,n}^{j+1}}\sum_{i=1}^n
K_j \left(\frac{t-t_i}{\hat{\lambda}_{j,n}}\right)(t_i-t_{i-1})y_i,\ \ l=0,...,r,
\end{equation}
and plug \fr{eq:deriv_est} into \fr{eq:estimator} or \fr{eq:f_explicit_est}.

Note that the resulting estimators $\hat{f}_n$ are inherently adaptive to the
smoothness of the underlying function $f$ in \fr{eq:model} which is rarely
known in practice.


\subsection{Adaptive minimaxity}
\label{sec:adaptive_minimax}

The following theorem establishes the upper bound for the $L_2([0,T_n])$-risk
of the estimator $\hat{f}_n$ defined in Section \ref{subsec:noisy_solution}
over Sobolev classes:

\begin{theorem} \label{th:risk_upper_bounds}
Let condition \fr{k_cond} and Assumptions (A1)-(A4) hold. Consider kernels
$K_j,\;j=0,\ldots,r$ of orders $(L,j),\;L>r$ satisfying the conditions (K1) and (K2).
Let $\hat{f}_n$ be the estimator of $f$ of the form \fr{eq:estimator} or
\fr{eq:f_explicit_est}, where $\widehat{q^{(j)}}(t)$'s are given by
\fr{eq:deriv_est}.
Then, for all $1 \leq m \leq \min(L,\nu+1)-r$, and $A>0$, one has
\be \label{eq:upper_finite}
\sup_{f \in W^m(A)} E\|\hat{f}_n-f\|^2_{L_2([0,T_n])}=O\left(\left(\frac{T_n^2}{n}\right)^{\frac{2m}{2(m+r)+1}}\right).
\ee
\end{theorem}
Under the additional conditions on $f$ and $T_n$, the results of Theorem
\ref{th:risk_upper_bounds} can be easily extended to the entire nonnegative
half-line:

\begin{corollary} \label{cor:risk_halfline}
Let conditions of Theorem \ref{th:risk_upper_bounds} hold and also
there exists $\rho \geq 1$ such that
$\int_0^\infty  t^{2 \rho} f^2 (t) dt < \infty$
and
$\lim_{n \rightarrow \infty} T_n^{-2\rho} n  < \infty$.
Let $\hat{f}_n$ be as in Theorem \ref{th:risk_upper_bounds}  for $t \leq T_n$ and
$\hat{f}_n \equiv 0$ for $t>T_n$. Then,
$$
\sup_{f \in W^m(A)} E\|\hat{f}_n-f\|^2_{L_2([0,\infty))}=O\left(\left(\frac{T_n^2}{n}\right)^{\frac{2m}{2(m+r)+1}}\right)
$$
for all  $1 \leq m \leq \min(L,\nu+1)-r$ and $A>0$.
\end{corollary}

Note that the upper bounds established in Theorem \ref{th:risk_upper_bounds} and
Corollary \ref{cor:risk_halfline} coincide with the minimax lower bound for the risk
obtained in Theorem \ref{th:lower_bounds} and, thus, cannot be improved.
Hence, the derived Laplace deconvolution estimators are asymptotically adaptively
minimax over entire range of Sobolev classes.


\section{Examples and simulation study }
\label{sec:example_simul}
\setcounter{equation}{0}

\subsection{Examples of explicit Laplace deconvolution estimators}
\label{sec:examples}

In what follows, we shall consider two examples of  construction of explicit
estimators of $f$ in the Laplace convolution problem.\\

{\bf Example 1.} Consider \fr{eq:model} with
$$
g(t) = (bt - \sin(bt)) e^{-at},\ \ a>0.
$$
It is easy to see that $r=4$ and $B_4 = b^3$ in \fr{k_cond}, and $\tilde{g}$ is of the form
\be \label{eq:tilG1}
\tilde{g} (s) = b^3 (s+a)^{-2} ((s+a)^2 + b^2)^{-1}.
\ee
Hence,  $\tilde{g} (s)$ has no zeros and one can use Theorem \ref{th:semiexplicit}
for recovering and estimating $f$. By \fr{eq:Lapl_phi} one has
\begin{equation*}
\widetilde{\phi} (s) = - \left(\frac{4a}{s} + \frac{6 a^2 + b^2}{s^2}
+ \frac{4a^3 + 2ab^2}{s^3} + \frac{a^4 + a^2b^2}{s^4} \right),
\end{equation*}
so that, in  \fr{eq:f_explicit} and \fr{eq:f_explicit_est}, one has
$\alpha_{0,0} = -4a$, $\alpha_{0,1} = -(6 a^2 + b^2)$,
$\alpha_{0,2} = -(4a^3 + 2ab^2)$, $\alpha_{0,3} = -(a^4 + a^2b^2)$ and $\phi_1(x) =0$.
Hence, using \fr{eq:f_explicit_est}
$$
\hat{f}_n(t)  =  b^{-3} \lkv \widehat{q^{(4)}}_{\hat{\lambda}_{n,4}}  (t) +
4a \widehat{q'''}_{\hat{\lambda}_{n,3}}(t) + (6 a^2 + b^2) \widehat{q''}_{\hat{\lambda}_{n,2}} (t)
+ (4a^3 + 2ab^2) \widehat{q'}_{\hat{\lambda}_{n,1}} (t) + (a^4 + a^2b^2) \widehat{q}_{\hat{\lambda}_{n,0}} (t) \rkv,
$$
where  $\hat{\lambda}_{n,l}$, $l=0,1,..., 4$, are defined in \fr{eq:adapband}.
The rate of convergence of $\hat{f}_n$ over $W^m$ is given by
\fr{eq:upper_finite} with $r=4$ and is $O\left(\left(\frac{T_n^2}{n}\right)^{\frac{2m}{2m+9}}\right)$.

{\bf Example 2.} Consider \fr{eq:model} with
\begin{equation} \label{eq:AIF_kerenel}
g(t) =   e^{-at} t^{r-1} \sum_{j=0}^k \frac{\rho_j}{(j+r-1)!} t^j,\ \ a>0,
\end{equation}
where $k \geq 0$ and $r \geq 1$ are integers and $\rho_0=1$. In this case,
\fr{k_cond} holds with  $B_r =1$ and
$$
\tilde{g}(s) = (s+a)^{-(k+r)}{\cal P}(s),
$$
where 
\be \label{eq:Ps} 
{\cal P}(s)=\sum_{j=0}^k \rho_j (s+a)^{k-j}.
\ee
Therefore,
$$
\widetilde{\phi}(s)=\frac{s^r {\cal P}(s)-(s+a)^{k+r}}{s^r {\cal P}(s)},
$$
In particular, for $k=0$ and $r=1$, ${\cal P}(s)$ has no roots, so that
$b_0 = -a$ and  we recover the result of
Dey, Martin and Ruymgaart (1998): $f(x) = q'(t) + a q(t)$.
For $k=1$, $\rho_0=1$ and $\rho_1=b$, one has $g(t) =  e^{-at} (bt+1)$ and
$\tilde{g}(s) = (s+a)^{-2} (s+a+b)$,
so that ${\cal P}(s)$ has a single root $s_1 = -(a+b)$ of multiplicity $\alpha_1 =1$.
Hence, $b_0 = b^2(a+b)$, $a_{1,0} = -b^2$ and $\phi_1(x) = -b^2 (a+b)^{-1} e^{-(a+b)x}$ in formula \fr{eq:f_explicit} leading
to the estimator of $f$ of the form
\be \label{eq:est_simul_gen}
\hat{f}_n(t)  =  \widehat{q'}_{\hat{\lambda}_{n,1}} (t)  + (a-b) \widehat{q}_{\hat{\lambda}_{n,0}} (t)
+ b^2  \int_0^t \widehat{q}_{\hat{\lambda}_{n,0}} (t-x) e^{-(a+b)x} dx.
\ee
The asymptotic minimax rate of convergence of $\hat{f}_n$ in \fr{eq:est_simul_gen}
over $W^m$ is $O\left(\left(n^{-1}\, T_n^2 \right)^{\frac{2m}{2m+3}}\right)$.

For   general values of  $k$ and $r$, the exact form of the solution \fr{eq:f_explicit} strongly depends on the roots
of the polynomial ${\cal P} (s)$ given by \fr{eq:Ps}. 
Assume that ${\cal P}(s)$ has $k$ distinct roots. 
\ignore{
$s_{2l-1} = u_{l1} + i u_{l2}$ and $s_{2l}  = u_{l1} - i u_{l2}$,
$l=1, \ldots, m$, and $k-2m$ real roots $s_l, \ldots, s_{k}$, $l=k-2m+1, 
\ldots, k$ and that all $k$ roots are distinct. 
}
Then,
${\cal P}(s) = \prod_{l=1}^k   (s-s_l)$  and  $1/\tilde{g} (s)$ allows a partial fraction decomposition
\begin{equation} \label{eq:partial_frac}
\frac{1}{\tilde{g} (s)} = \sum_{j=0}^r \alpha_j (s+a)^j + \sum_{l=1}^{k} \frac{\beta_l}{s - s_l}.
\end{equation}
By observing that $\sum_{j=0}^r \alpha_j s^j$ is the quotient of $s^{r+k}$ and $ \sum_{j=0}^k \rho_j s^{k-j}$,
one can recursively evaluate $\alpha_j$, $j=1, \ldots, r$, in \fr{eq:partial_frac} as
$$
\alpha_r =1,\ \ \ \alpha_{r-l} = - \sum_{j =\max(0, l-k)}^{l-1} \alpha_{r-j} \rho   _{l-j},\ \ l=1, \ldots, r.
$$
The values of $\beta_l$ can be obtained by multiplying both sides of equation 
\fr{eq:partial_frac} by ${\cal P} (s) /(s - s_l)$
and  setting $s=s_l$:
$$
\beta_l = (s_l + a)^{k+r} \prod_{\stackrel{j=1}{j \neq l}}^{k} (s_l-s_j)^{-1},\ \ l=1,\ldots,k.
$$
The respective expression for $f$ is of the form $f=f_1 + f_2$ where
\begin{eqnarray} \label{eq:ex_explicit_solution}
f_1(t) & = & q^{(r)} (t) + \sum_{l=0}^{r-1} q^{(l)} (t) \sum_{j=l}^r {j \choose l} a^{j-l} \alpha_j,\\
f_2(t) & = & \sum_{l=1}^k \beta_l \int_0^t e^{s_l x} q(t-x) dx,
\end{eqnarray}
which can easily be reduced to representation \fr{eq:f_explicit}.
\ignore{
Note that $s_l$ and $\beta_l$, $l = 2m+1, \ldots, k$,  have real values while $s_{2l-1}$ and $s_{2l}$ as well as
$\beta_{2l-1}$ and $\beta_{2l}$, $l=1, \ldots, m$, are complex conjugates of each other. Hence, one can rewrite expression for $f_2$ as
$$
f_2 (x) = \sum_{l=1}^m   \int_0^t e^{u_{l1} x} [\gamma_{l1} \cos(u_{l2} x) + \gamma_{l2} \sin(u_{l2} x)] q(t-x) dx
+    \sum_{l=2m+1}^k \beta_l \int_0^t e^{s_l x} q(t-x) dx
$$
where $\gamma_{l1} = 2 \mbox{Re} (\beta_{2l})$ and $\gamma_{l2} = 2 \mbox{Im} (\beta_{2l})$.
}

Under assumptions (A2) and (A3), the asymptotic minimax rate of convergence of the estimator \fr{eq:est_simul_gen}
is provided by Theorem \ref{th:risk_upper_bounds}.


\subsection{Simulation study}
\label{sec:simul}

In this section we present the results of a simulation study to illustrate
finite sample performance of the Laplace deconvolution procedure developed
above.

\begin{figure}[ht]
\[\includegraphics[height=3cm]{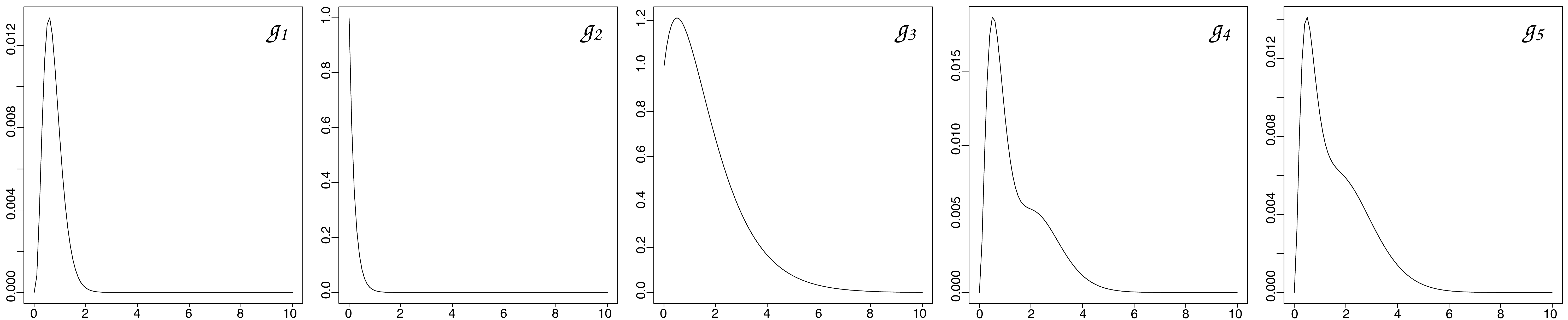} \]\vspace{-10mm}
\caption{\small From left to right: the kernels $g_1$ to $g_5$. \label{fig:g-functions}}
\end{figure}

First, we consider the data simulated according to the model (\ref{eq:model}) with
five convolution kernels $g_1$, \ldots, $g_5$, where
$$
g_1(t)=e^{-5t} (2t-\sin(2t)),\qquad g_2(t)=e^{-5t},\qquad g_3(t)=e^{-t} (2t+1).
$$
Kernel $g_1$ mimics an ideal behaviour of AIF in the DCE-imaging
(see Example 1 in Section \ref{sec:form_motive}), while $g_2$ and
$g_3$ are examples of kernels considered in Example 2 from Section
\ref{sec:examples}. In particular, $g_2$ corresponds to Dey, Martin
and Ruymgaart (1998) framework. Kernels $g_4$ and $g_5$ also fall
within the general form of Example 2 from Section \ref{sec:examples}
with $r=3$ and were defined by the $k$ roots $s_1$, \ldots, $s_k$ of
the polynomial ${\cal P}(s)$ in \fr{eq:Ps}. For $g_4$ we considered
four roots $(-4 \pm 2.5i,  -0.75\pm 1.5i)$, while for $g_5$ we
added two more conjugate roots $ -2 \pm 2i$. Both $g_4$ and $g_5$
can be seen as more realistic scenarios in the DCE-imaging. All the
five kernels are presented on Figure \ref{fig:g-functions}.

\begin{figure}[ht]
\[\includegraphics[height=3cm]{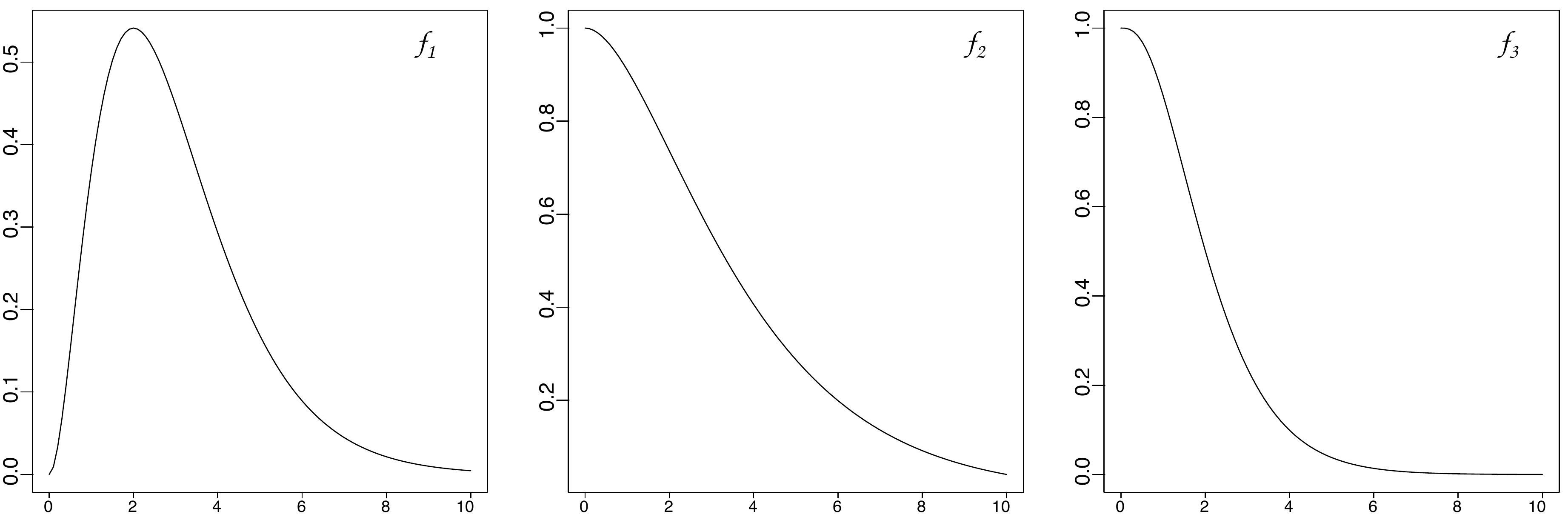} \]\vspace{-10mm}
\caption{\small The true unknown function   $f$ \label{fig:f-functions} from left to right: $f_1$ to $f_3$.}
\end{figure}

The chosen true functions $f$ in \fr{eq:model} are
$f_1(t)=t^2 e^{-t}$, $f_2(t) = 1-\Gamma_{2,2}(t)$ and $f_3(t) = 1-\Gamma_{3,0.75}(t)$, where
$\Gamma_{\alpha,\theta}$ is the c.d.f of the Gamma
distribution with the shape parameter $\alpha$ and the scale parameter $\theta$
(see Figure \ref{fig:f-functions}).
Functions $f_2$ and $f_3$ mimic sojourn time distributions of the particles
of a contrast agent in DCE-imaging experiments,
while $f_1$ is aimed to be a more general case.

The Laplace convolution $q=f \ast g$ which produces observations  in
(\ref{eq:model}) has been numerically computed using trapezoidal rule for
approximation of the integral. The noise levels for each of the kernels $g_1$, ..., $g_5$
was chosen as  $\sigma_0(g_j)/2^i$, $i=0,\ldots 4;\; j=1,\ldots,5$,
where the nominal noise levels $\sigma_0(g_j)$ were
0.001, 0.1, 0.01, 0.002, 0.002 for $g_1, \ldots, g_5$ respectively.
We ran  simulations with $n=100$ and $n=250$ and  regular design for the $t_i$ equally spaced  between 0 and
$T_n=10$.


Following  construction in Gasser, M\"uller and
Mammitzsch (1985), we derived kernels $K_j$ of orders $(L, j)$ for estimating
the derivatives $q^{(j)}$ of $q$, $j=0, \ldots, r$ for various values of $L$.
In our simulations we used $L=8$ as an upper bound of the regularity of the
kernel  since  higher values of $L$ lead to  numerically unstable computations
and/or provide very little advantage in terms of precision. Finally, we used
boundary kernels in order to stabilized the computations as suggested in
Gasser, M\"uller and Mammitzsch (1985). In all simulations, due to the regular 
fixed design, 
$\mu=1$ in Assumption (A4). We chose $a=1.2$ in \fr{LambdaDef} and 
$C_j=1$ in \fr{c0}. 
Since the constant $4$ in the Lepski's threshold in \fr{eq:adapband} is known to 
be too large for practical applications, we tried several values and 
``tuned'' it to 3. 

\begin{figure}[ht]
\[\includegraphics[width=0.99\textwidth]{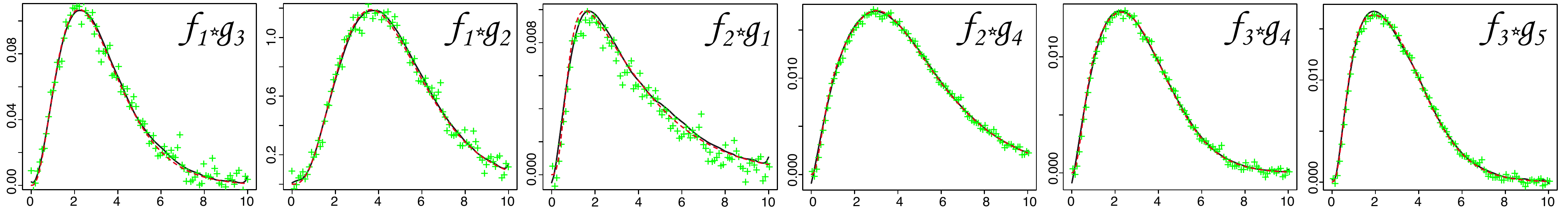} \]\vspace{-10mm}
\caption{\small Laplace convolution $q$ (dotted red  line) of the known kernel $g$ and the unknown function   $f$,
$n=100$ noisy observations of $q$ (green pluses) and estimated value of $q$ (black line). The choice of $f$ and $g$ used for each simulation are specified on each sub-figure by the convolution product. The noise level has been chosen as
$\sigma_0(g_j)/2$ for the three left figures and  $\sigma_0(g_j)/8$ for the three right figures.
\label{fig:q-estimation}}
\end{figure}

\begin{figure}[ht]
\[\includegraphics[width=0.99\textwidth]{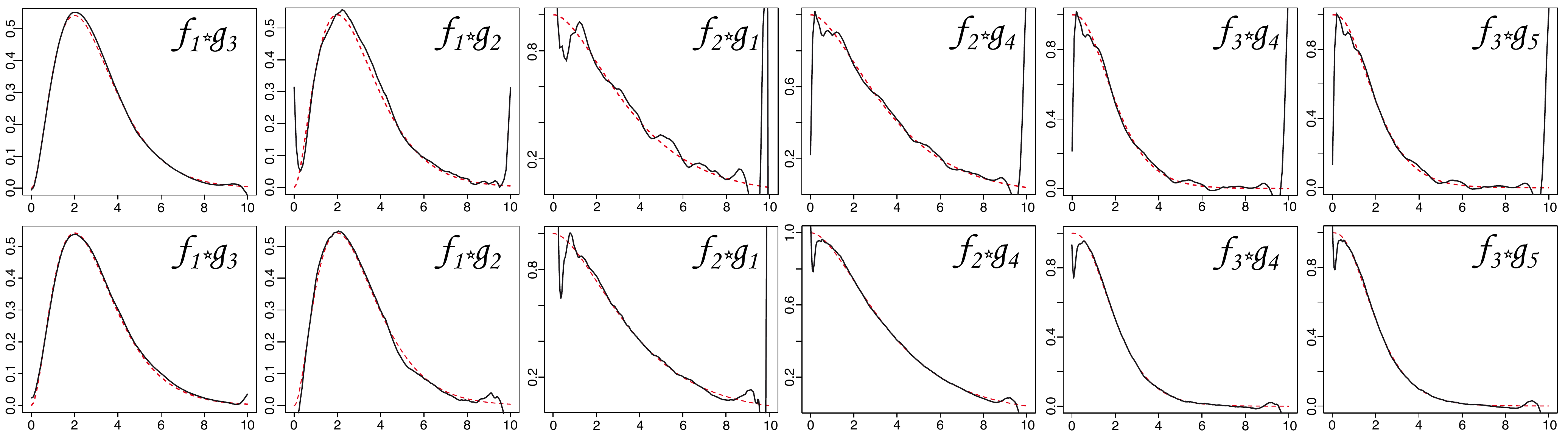} \]\vspace{-10mm}
\caption{\small True unknown $f$ (red dotted line) and its estimate (plain black line) for $n=100$ (top line) and $n=250$ (bottom line).
\label{fig:f-estimation}}
\end{figure}

\begin{figure}[ht]
\[\includegraphics[height=10cm,width=15cm]{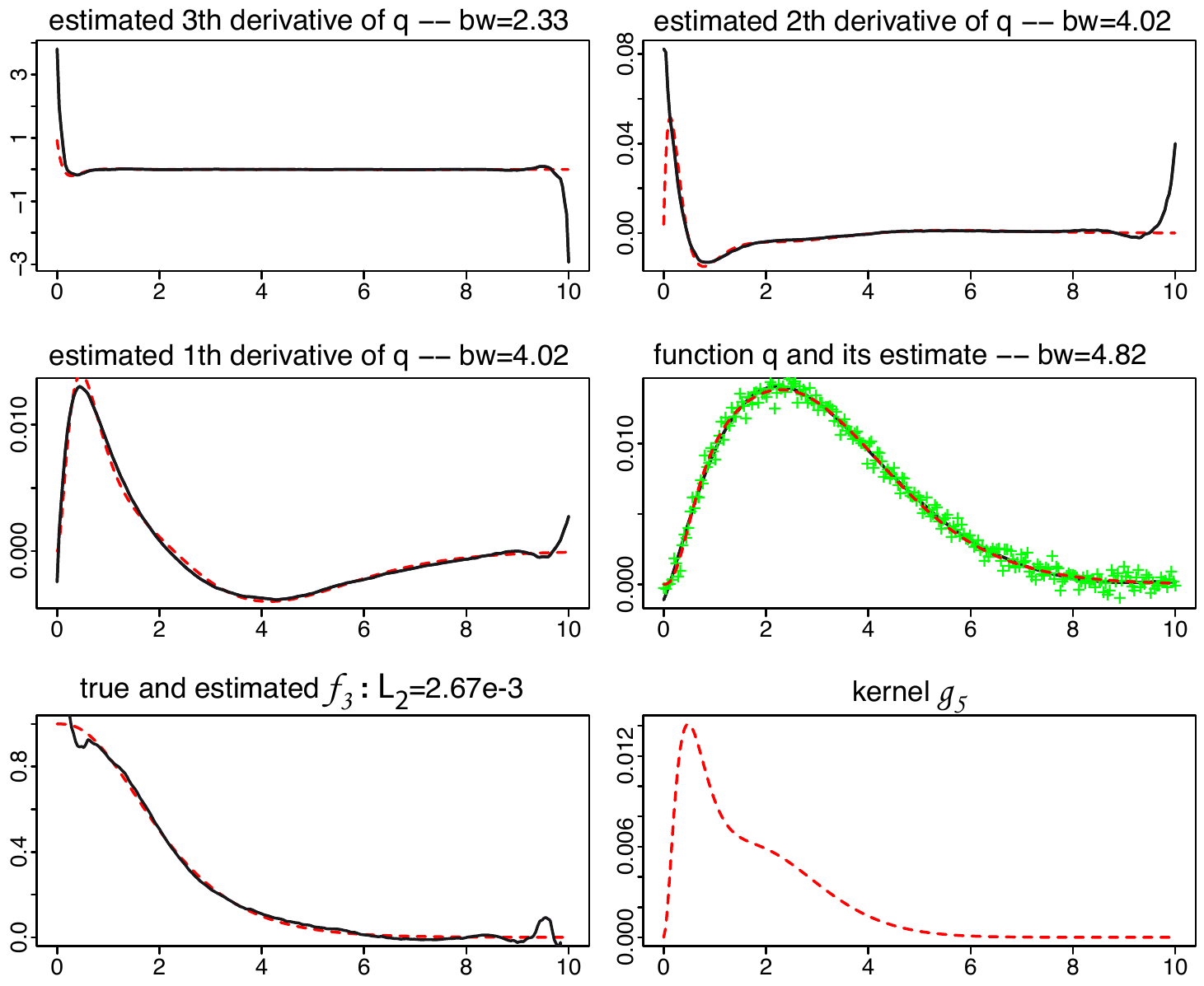} \]\vspace{-10mm}
\caption{\small Example of estimation  of unknown function $f_3$ with kernel $g_5$ using $n=250$ and $\sigma_0(g_5)/4$.
Here $r=3$. The top four sub-figures show true   $q$ and its three first derivatives $q^{(s)}$, $s=1,2,3$ (dotted red lines)
and their estimators  (plain black lines). Selected bandwidths are specified for each estimator. Bottom left: the true function $f$ (dotted red line) and its estimate (plain black line). Bottom right: the kernel $g_5$.
\label{fig:realistic}}
\end{figure}

Figures \ref{fig:q-estimation} and \ref{fig:f-estimation} provide  examples of deconvolution estimators based
on single samples.
Figure \ref{fig:realistic} illustrates that deconvolution 
estimators show good precision although boundary effects in estimating 
high-order derivatives remain despite  the use of boundary kernels.

\begin{figure}[ht]
\[\includegraphics[height=10cm]{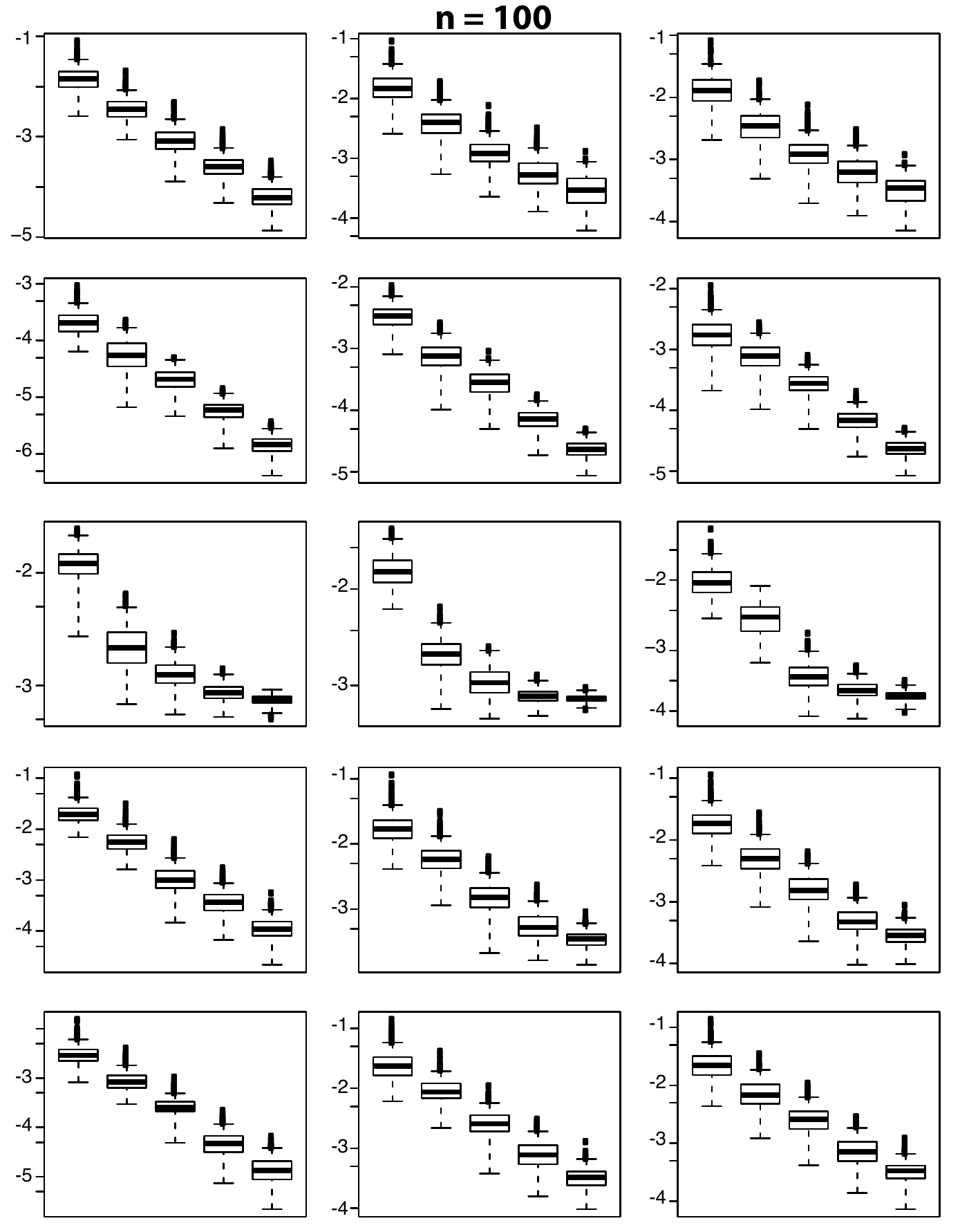}\qquad\includegraphics[height=10cm]{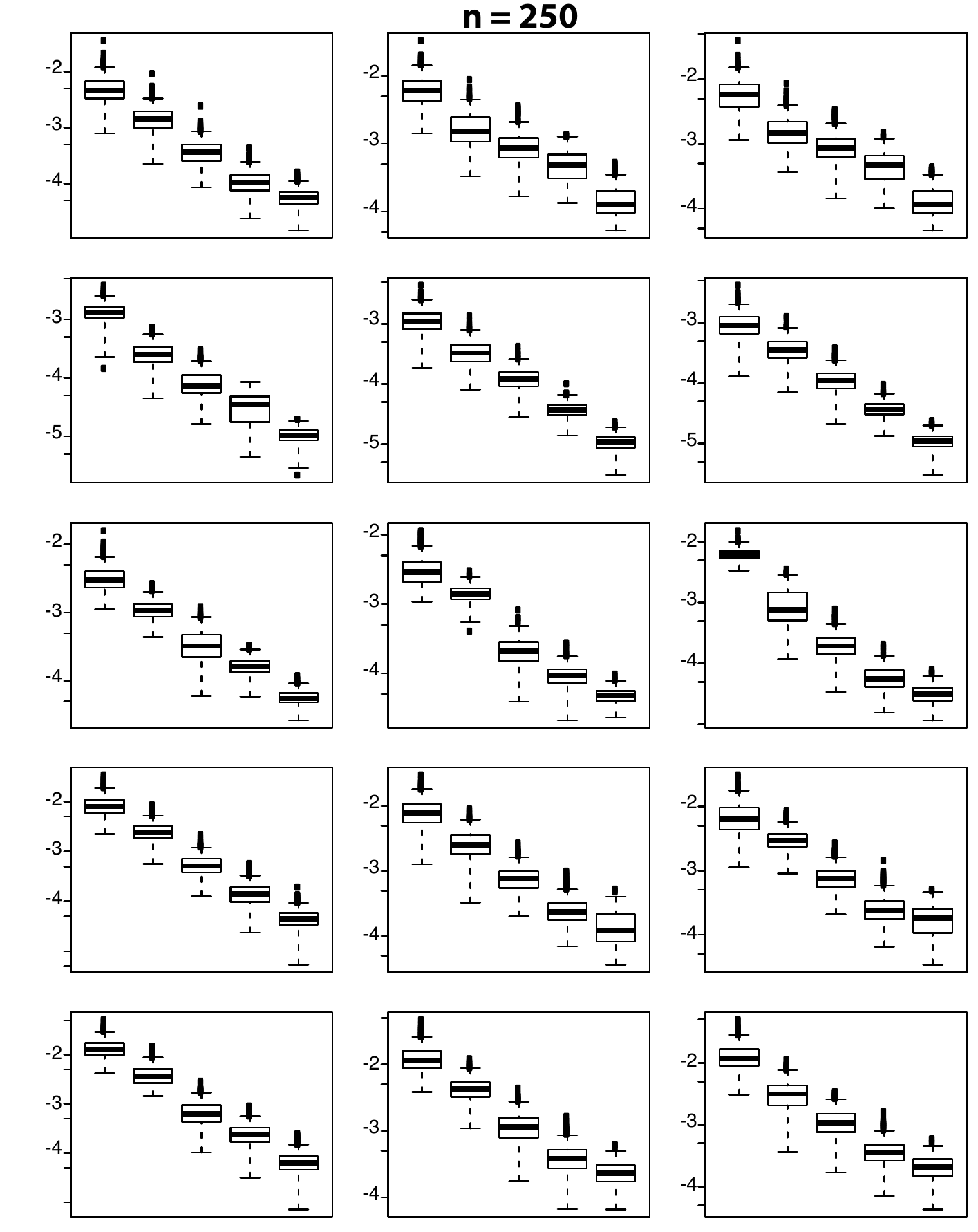} \]\vspace{-10mm}
\caption{ Box-plots of the mean square errors for $n=100$ (left) and $n=250$ (right), 400 simulation runs
 and each of the triplets $(g,f,\sigma_0(g))$. Each line represents a kernel from $g_1$ (top) to $g_5$ (bottom).
Each column represents an unknown function $f$ from $f_1$ (left) to $f_3$ (right). In every sub-figure,
going from left to right, each boxplot corresponds to a different noise level $\sigma_0(g_j)/2^i$ for $i=0,\ldots,4$.
The empirical risks are presented on log-scale with the basis 10.
\label{fig:boxplots}}
\end{figure}

For each combination of true function $f$, kernel $g$, sample size
$n$ and the noise level, we ran 400 simulations and calculated mean
square errors. In order to remove the influence of boundary effects
(see comments above), we did not include 20\% of the boundary points (10\% at
each boundary). The box-plots of the resulting mean square
errors are presented on Figure \ref{fig:boxplots}. Table
\ref{tab:risks} shows the average mean square errors and standard
deviations (in parentheses) over 400 simulation runs.

\begin{table}\centering\small
\begin{tabular}{llccccc|c}
\multicolumn{2}{l}{$\hat R(\hat f)$}             & i=0 & i=1 & i=2 & i=3 & i=4 \\\hline
\multirow{3}{*}{$g_1$} & $f_1$ & 1.6e-2 (9.5e-3) & 4.1e-3 (2.5e-3) & 9.9e-4 (6.6e-4) & 2.9e-4 (1.7e-4) & 7.6e-5 (4.7e-5) &\multirow{15}{*}{\begin{rotate}{-90}\boldmath$n=100$\end{rotate}}\\
 & $f_2$ & 1.7e-2 (1.0e-2) & 4.5e-3 (2.8e-3) & 1.4e-3 (8.4e-4) & 6.7e-4 (4.2e-4) & 3.4e-4 (1.9e-4) \\
 & $f_3$ & 1.5e-2 (9.6e-3) & 4.0e-3 (2.6e-3) & 1.4e-3 (9.0e-4) & 7.4e-4 (4.4e-4) & 3.6e-4 (1.7e-4) \\\cline{1-7}
\multirow{3}{*}{$g_2$} & $f_1$ & 2.3e-3 (1.1e-3) & 6.9e-4 (4.6e-4) & 2.2e-4 (8.8e-5) & 6.1e-5 (2.2e-5) & 1.5e-5 (5.5e-6) \\
 & $f_2$ & 3.5e-3 (1.5e-3) & 8.2e-4 (4.0e-4) & 2.9e-4 (1.3e-4) & 7.5e-5 (2.8e-5) & 2.4e-5 (7.7e-6) \\
 & $f_3$ & 2.1e-3 (1.4e-3) & 8.6e-4 (4.4e-4) & 2.9e-4 (1.1e-4) & 7.2e-5 (2.8e-5) & 2.5e-5 (7.8e-6) \\\cline{1-7}
\multirow{3}{*}{$g_3$} & $f_1$ & 1.2e-2 (3.7e-3) & 2.4e-3 (1.0e-3) & 1.3e-3 (3.6e-4) & 8.7e-4 (1.4e-4) & 7.4e-4 (7.2e-5) \\
 & $f_2$ & 1.4e-2 (3.8e-3) & 3.5e-3 (9.2e-4) & 2.2e-3 (5.0e-4) & 1.7e-3 (1.9e-4) & 1.6e-3 (9.1e-5) \\
 & $f_3$ & 1.0e-2 (3.6e-3) & 4.3e-3 (1.4e-3) & 1.2e-3 (3.7e-4) & 8.2e-4 (1.5e-4) & 7.1e-4 (7.3e-5) \\\cline{1-7}
\multirow{3}{*}{$g_4$} & $f_1$ & 2.2e-2 (1.2e-2) & 6.3e-3 (3.3e-3) & 1.2e-3 (8.3e-4) & 4.2e-4 (2.4e-4) & 1.2e-4 (6.5e-5) \\
 & $f_2$ & 2.0e-2 (1.2e-2) & 6.6e-3 (3.6e-3) & 1.8e-3 (1.0e-3) & 6.0e-4 (3.1e-4) & 3.6e-4 (1.1e-4) \\
 & $f_3$ & 2.1e-2 (1.2e-2) & 5.7e-3 (3.5e-3) & 1.8e-3 (1.0e-3) & 5.5e-4 (3.0e-4) & 3.0e-4 (1.1e-4) \\\cline{1-7}\multirow{3}{*}{$g_5$} & $f_1$ & 3.2e-2 (1.7e-2) & 9.4e-3 (4.6e-3) & 2.8e-3 (1.2e-3) & 5.3e-4 (3.2e-4) & 1.6e-4 (9.9e-5) \\
 & $f_2$ & 2.7e-2 (1.7e-2) & 9.8e-3 (4.5e-3) & 2.9e-3 (1.5e-3) & 8.8e-4 (4.7e-4) & 3.4e-4 (1.4e-4) \\
 & $f_3$ & 2.6e-2 (1.7e-2) & 8.1e-3 (4.8e-3) & 2.9e-3 (1.5e-3) & 8.3e-4 (4.8e-4) & 3.4e-4 (1.4e-4) \\\hline\hline
\multirow{3}{*}{$g_1$} & $f_1$ & 5.5e-3 (3.4e-3) & 1.6e-3 (9.0e-4) & 4.1e-4 (2.3e-4) & 1.2e-4 (5.7e-5) & 5.9e-5 (2.2e-5) & \multirow{15}{*}{\begin{rotate}{-90}\boldmath$n=250$\end{rotate}}\\
 & $f_2$ & 7.0e-3 (3.8e-3) & 1.9e-3 (1.2e-3) & 1.0e-3 (5.6e-4) & 5.3e-4 (2.6e-4) & 1.6e-4 (9.1e-5) \\
 & $f_3$ & 6.5e-3 (3.9e-3) & 1.8e-3 (1.1e-3) & 1.0e-3 (5.4e-4) & 5.0e-4 (2.6e-4) & 1.5e-4 (8.4e-5) \\\cline{1-7}
\multirow{3}{*}{$g_2$} & $f_1$ & 1.4e-3 (5.2e-4) & 2.7e-4 (1.2e-4) & 9.1e-5 (5.4e-5) & 3.5e-5 (1.9e-5) & 1.0e-5 (3.1e-6) \\
 & $f_2$ & 1.2e-3 (5.4e-4) & 3.7e-4 (1.8e-4) & 1.3e-4 (5.2e-5) & 3.8e-5 (1.1e-5) & 1.1e-5 (3.2e-6) \\
 & $f_3$ & 1.0e-3 (4.9e-4) & 3.9e-4 (1.8e-4) & 1.2e-4 (4.8e-5) & 3.8e-5 (1.1e-5) & 1.1e-5 (3.4e-6) \\\cline{1-7}
 \multirow{3}{*}{$g_3$} & $f_1$ & 3.4e-3 (1.7e-3) & 1.2e-3 (4.1e-4) & 3.8e-4 (2.2e-4) & 1.7e-4 (4.7e-5) & 5.8e-5 (1.4e-5) \\
 & $f_2$ & 3.6e-3 (2.2e-3) & 1.4e-3 (3.8e-4) & 2.3e-4 (1.1e-4) & 9.6e-5 (3.4e-5) & 4.8e-5 (1.2e-5) \\
 & $f_3$ & 6.3e-3 (1.4e-3) & 1.0e-3 (7.1e-4) & 2.1e-4 (1.1e-4) & 6.1e-5 (2.9e-5) & 3.3e-5 (1.2e-5) \\\cline{1-7}
\multirow{3}{*}{$g_4$} & $f_1$ & 8.9e-3 (4.5e-3) & 2.7e-3 (1.1e-3) & 5.8e-4 (3.0e-4) & 1.6e-4 (8.6e-5) & 4.9e-5 (2.3e-5) \\
 & $f_2$ & 8.8e-3 (4.6e-3) & 2.9e-3 (1.5e-3) & 8.2e-4 (3.8e-4) & 2.8e-4 (1.5e-4) & 1.6e-4 (9.6e-5) \\
 & $f_3$ & 7.7e-3 (4.8e-3) & 3.2e-3 (1.3e-3) & 8.1e-4 (3.7e-4) & 2.9e-4 (1.9e-4) & 1.9e-4 (9.6e-5) \\\cline{1-7}
\multirow{3}{*}{$g_5$} & $f_1$ & 1.5e-2 (7.0e-3) & 4.1e-3 (2.0e-3) & 7.4e-4 (4.4e-4) & 2.7e-4 (1.4e-4) & 7.2e-5 (4.0e-5) \\
 & $f_2$ & 1.3e-2 (6.7e-3) & 4.6e-3 (1.9e-3) & 1.3e-3 (6.4e-4) & 4.3e-4 (2.3e-4) & 2.4e-4 (9.9e-5) \\
 & $f_3$ & 1.4e-2 (6.9e-3) & 3.5e-3 (1.9e-3) & 1.2e-3 (5.9e-4) & 4.1e-4 (2.2e-4) & 2.2e-4 (9.4e-5) \\\hline
\end{tabular}
\caption{\small Average (over 400 simulation runs) mean square errors and standard deviations (in parentheses)
for  kernels $g_1$ to $g_5$ and unknown function $f_1$ to $f_3$, for $n=100$ (upper part) and $n=250$ (lower part) and for the noise level equal to
$\sigma_0(g_j)/2^i$, $j=1,\ldots,5;\;i=0,\ldots,4$.
\label{tab:risks}}
\end{table}

\section{Discussion}
\label{sec:discussion}
In the present paper, we consider Laplace deconvolution problem with discrete
noisy data observed on the interval whose length $T_n$ may increase with the sample
size $n$. Although this problem arises in a variety of applications, to the best of our
knowledge, it has been given very little attention by the statistical community.
Our objective was to fill this gap and to  provide   statistical analysis
of Laplace deconvolution problem with noisy discrete data.

The main contribution of the paper is explicit construction of
a rate-optimal (in the minimax sense) Laplace deconvolution
estimator which is adaptive to the regularity of the unknown function.
We show that the original Laplace deconvolution problem can be  reduced to
nonparametric  estimation of a regression function and its derivatives on the interval of
growing length $T_n$. Although the latter problem has been well studied on a finite interval, the
asymptotic increase of its length  as the sample size grows raises a
new challenge. Whereas the forms of the estimators   remains standard, the choices of
the parameters and the minimax convergence rates, which are expressed in terms of $T_n^2/n$
in this case,  are affected by the asymptotic growth of the length of the interval.

In the present paper, we use kernel estimators with a global bandwidth adaptively chosen by
the Lepski procedure (e.g., Lepski, 1991) and establish  asymptotic minimaxity
of the resulting Laplace deconvolution estimator over a wide range of
Sobolev classes. One can, however, apply  other types of estimators
(e.g., local polynomial regression, splines or wavelets).
In particular, we believe that the use of wavelet-based methods can extend the
adaptive minimaxity range from Sobolev to more general Besov classes.

We illustrate the theory by examples of construction of explicit expressions for
estimators of $f$ based on  observations governed by equation \fr{eq:model}
with various kernels. Simulation study shows, that, in addition to providing
asymptotic optimality, the proposed Laplace deconvolution estimator
demonstrates   good  finite sample performance.

The present paper  provides the first comprehensive statistical
treatment of Laplace deconvolution problem, though
a number of open questions remain beyond its scope.
In particular, an interesting challenge would be to study Laplace
deconvolution with an unstable resolvent, where Assumption (A2) does not hold.
Another important problem would be to study the equation \fr{eq:model}
when the kernel $g$ is not completely known and is estimated from observations.

\section*{Acknowledgments}

Marianna Pensky was partially supported by National Science Foundation
(NSF), grant  DMS-1106564. We would like to thank Alexander Goldenshluger and Oleg Lepski for fruitful discussions
of the paper.


\section{Appendix}
\label{sec:append}
\setcounter{equation}{0}

Throughout the proofs we use $C$ to denote a generic positive constant, not
necessarily the same each time it is used, even within a single equation.

\subsection*{Proof of Theorem \ref{th:lower_bounds}}

 Although the rates are derived by standard methods described in, e.g., Tsybakov (2009),
the challenging part of the proof is constructing the set of test functions and, subsequently,
  producing upper bounds for the Kullback-Leibler divergence.

The main idea of the proof is to find a subset of functions ${\cal F}\subset W^{m}(A)$
 such that  for any pair $f_1,\;f_2 \in {\cal F}$,
\be \label{eq:cond1}
\|f_1-f_2\|^2_{L_2([0,T_n])} \geq 4 C (T^2_n n^{-1})^{2m/(2(m+r)+1)}
\ee
and the Kullback-Leibler divergence
\be \label{eq:cond2}
\KK (\PP_{f_1}, \PP_{f_2})=
\frac{||\bq_1-\bq_2||^2_{\mathbb{R}^n}}{2\sigma^2} \leq
\frac{\log  {\rm card}({\cal F})}{16},
\ee
where $\log$ stands for natural logarithm and vectors $\bq_j$,  $j=1,2$, have components
$\bq_{ji}=(g \ast f_j)(t_i)$, $i=1,...,n$.
The result will then follow immediately from
Lemma A.1 of Bunea, Tsybakov and Wegkamp (2007):

\begin{lemma} \label{Tsybakov} [Bunea, Tsybakov, Wegkamp~(2007), Lemma A.1]
Let   ${\cal F}$ be a set of functions of cardinality $\card({\cal F}) \geq 2$ such that \\
(i) $\|f_1-f_2\|^2 \geq 4 \delta^2 \quad$ for any $f_1, f_2 \in {\cal F}$, \ $f_1 \neq f_2$,\\
(ii) the Kullback divergences $\KK (\PP_{f_1}, \PP_{f_2})$ between the measures $\PP_{f_1}$ and
$\PP_{f_2}$ satisfy the inequality  $\KK (\PP_{f_1}, \PP_{f_2}) \leq  (1/16)\log(\card({\cal F})) \quad$
for any $f_1, f_2 \in {\cal F}$.\\
Then, for some absolute positive constant $C$,
 \begin{eqnarray*}
 \inf_{\tilde{f}_n} \sup_{f \in \ {\cal F}} \EE_f  \|\tilde{f}_n  -f\|^2  \geq C \delta^2,
 \end{eqnarray*}
where the infimum is taken over all estimates $\tilde{f}_n$ of $f$.
\end{lemma}

Without loss of generality, let us assume that the points are equally spaced,
i.e. $t_i - t_{i-1} = T_n/n$, $i=1,...,n$. To construct such a subset
${\cal F}$, define    integers  $M_n \geq 8$  and $N = \lkv \frac{n}{M_n} \rkv$,
the largest integer which does not exceed $n/M_n$. Let $\lamn = N T_n/n$ and define points
$z_l = l\, \lamn$, $l= 0, 1,..., M_n$. Note that the latter implies that
points of observation $t_j = j\, T_n/n$ in equation \fr{eq:model} are related to $z_l$
as $z_l = t_j$  where $j = N l$ for $l = 1,..., M_n$ and $j \leq N M_n$.
Note also that $\frac{T_n}{2 M_n} \leq \lamn \leq \frac{T_n}{M_n}$.

Let $k(\cdot)$ be an infinitely differentiable function with $supp(k)= [0,1]$
and such that
\be \label{k_conditions}
\int_0^1 x^j k(x) dx =0, j=0,..., r-1,\ \ \ \int_0^1 x^{r} k(x) dx \neq 0.
\ee
Introduce functions
$$
\varphi_j (x) = L\, \frac{\lamn^m}{\sqrt{T_n}}\,  k \lkr \frac{x - z_{j-1}}{\lamn} \rkr
\ \ \ l=1,..., M_n,
$$
where the constant $L>0$ will be defined later.
Note that $\varphi_j$ have non-overlapping supports, where
$supp(\varphi_j)= [z_{j-1}, z_j]$.

Consider the set of all binary sequences of the length $M_n \geq 8$:
$$
\Omega = \big\{ \bomega = (\omega_1, ..., \omega_{M_n}), \ \ \omega_j = \{0,1\} \big\} =  \{0,1\}^{M_n}
$$
and the corresponding subset of functions
\be \label{eq:setf}
{\cal F}=\{f_\omega : f_\omega(t)=\sum_{j=1}^{M_n} w_j \varphi_j(t),\;
\omega \in \tilde{\Omega}\}.
\ee
Here $\tilde{\Omega} \subset \Omega$ is such that
$\log_2 {\rm card}(\tilde{\Omega}) \geq  M_n/8  $ and the Hamming distance
$\rho(\bomega_1,\bomega_2)=\sum_{j=1}^{M_n} \mathbb{I}\{\bomega_{1j} \ne
\bomega_{2j}\} \geq M_n/8$ for any pair $\bomega_1,\bomega_2 \in \tilde{\Omega}$
(see, e.g., Lemma 2.9 of Tsybakov (2009) for construction of $\tilde{\Omega}$).

We now need to show that ${\cal F}$ in \fr{eq:setf} is exactly the required
set. Note first that since the supports of $\varphi_j$
are non-overlapping, for any $f_\omega \in {\cal F}$ a straightforward calculus yields
$$
||f_\omega||^2_{L_2([0,T_n])} \leq \sum_{j=1}^{M_n} ||\varphi_j||^2=
L^2 \frac{\lambda_n^{2s+1}}{T_n} M_n ||k||^2 = L^2 \lambda^{2m} ||k||^2 \leq
L^2 ||k||^2
$$
Similarly,
$$
||f^{(m)}_\omega||^2_{L_2([0,T_n])} \leq \sum_{j=1}^{M_n} ||\varphi^{(m)}_j||^2=
\frac{L^2}{T_n} m \lambda_n ||k^{(s)}||^2=L^2 ||k^{(m)}||^2 < \infty
$$
and therefore $f_\omega \in W^{(m)}(A)$, where $A=L||k||_{W^m}$.
Furthermore,
$$
||f_{\omega_1}-f_{\omega_2}||^2_{L_2([0,T_n])} =L^2 \frac{\lambda_n^{2m+1}}{T_n}
||k||^2 \rho(\bomega_1,\bomega_2) \geq L^2 \frac{\lambda_n^{2m+1}}{T_n}
\frac{M_n}{8} \geq 4 C \lambda_n^{2m}
$$
and \fr{eq:cond1} holds provided $\lambda_n \geq C  (T^2_n n^{-1})^{-1/(2(m+r)+1)}$
for some positive constant  $C$.

To verify \fr{eq:cond2}, note that
\be \label{eq:K_to_Q}
\KK (P_1, P_2) = \frac{1}{2\sigma^2} \  \sum_{i=1}^n \lkv q_1 (t_i) - q_2(t_i) \rkv^2 \leq  \frac{1}{ \sigma^2} \sum_{j=1}^2  Q(f_j)
\ee
where, suppressing index  $j$, we write
\begin{eqnarray*}
Q(f) &=&   \  \sum_{i=1}^n \lkv \int_0^{t_i} g(t_i - x) f  (x) dx \rkv^2 =
\frac{L^2 \lamn^{2m}}{T_n}\  \sum_{i=1}^n \lkv  \sum_{l=1}^{M_n} \omega_l^{(j)}
\int_0^{t_i} g(t_i - x) k \lkr \frac{x - z_{l-1}}{\lamn} \rkr dx \rkv^2.
\end{eqnarray*}

In order to obtain an upper bound for $Q(f)$ we need the following  supplementary lemma, the proof of which
is presented at the end of the section.

\begin{lemma}  \label{th:suppl_lemma}
Introduce functions $K_j(x)$ using the following recursive relation
\be \label{K_definition}
K_1(x) = \int_0^x k(t) dt,\ \ K_j (x) = \int_0^x K_{j-1} (t) dt,\ \ j=2, ..., r.
\ee
Then, under condition \fr{k_conditions}, functions $K_j(x)$,  $j=1,..., r$, are uniformly bounded and
$K_j (1) =0$, $j=1, \ldots, r$.
Moreover,
\begin{eqnarray}
\int_0^{t_i} g(t_i -x) k \lkr \frac{x-z_{l-1}}{\lamn} \rkr dx & = &
\lamn^r \lkv B_r K_r \lkr \frac{t_i - z_{l-1}}{\lamn} \rkr \II (z_{l-1} \leq y_i \leq z_l) \right. \nonumber\\
& + &
\left. \int_{\min(z_{l-1}, t_i)}^{\min(z_{l}, t_i)}  g^{(r)}(t_i -x) K_r \lkr \frac{x-z_{l-1}}{\lamn} \rkr dx \rkv.
\label{int_bound}
\end{eqnarray}
\end{lemma}

Applying equation \fr{int_bound} to the integral in $Q(f)$, obtain
\be \label{del12}
Q(f) \leq    2 L^2 \lamn^{2m + 2r} T_n^{-1} \, (\Delta_1 + \Delta_2)
\ee
where
\begin{eqnarray*}
\Delta_1 & = &  \sum_{i=1}^n \lkv  \sum_{l=1}^{M_n}
B_r\ K_r \lkr \frac{t_i - z_{l-1}}{\lamn} \rkr \II (z_{l-1} \leq y_i \leq z_l) \rkv^2,\\
\Delta_2 & = &  \sum_{i=1}^n \lkv  \sum_{l=1}^{M_n}
\int_{\min(z_{l-1}, t_i)}^{\min(z_{l}, t_i)}  g^{(r)}(t_i -x) K_r \lkr \frac{x-z_{l-1}}{\lamn} \rkr dx \rkv^2.
\end{eqnarray*}

Observe that for any $t$ and any $l_1$ and $l_2$ such that  $l_1 \neq l_2$, one has
$K_r (\lamn^{-1} (t-z_{l_1})) K_r (\lamn^{-1} (t-z_{l_2})) =0$. Also,
for each $i$, $K_r (\lamn^{-1} (ti -z_{l})) \neq 0$ for only one value of $l$, namely, for $l = [i/N] +1$
where $[x]$ is the largest integer which does not exceed $x$. Therefore,
\be \label{Delta1}
\Delta_1 \leq B_r^2  \sum_{i=1}^n  K_r^2 \lkr \frac{t_i - z_{[i/N]}}{\lamn} \rkr
\leq n\, \, B_r^2 \| K_r \|^2_\infty,
\ee
where $\| \cdot \|_\infty$ is the supremum norm. In order to obtain an upper bound for $\Delta_2$,
observe that for any nonnegative function $F(x)$ one has
$$
\int_{\min(z_{l-1}, t_i)}^{\min(z_{l}, t_i)}  F(x) dx \leq \int_{ z_{l-1}}^{ z_{l}}  F(x) dx.
$$
Hence,  we derive
\begin{eqnarray}
\Delta_2 & \leq &   \sum_{i=1}^n \lkv  \sum_{l=1}^{M_n} \int_{ z_{l-1}}^{ z_{l}} \left| g^{(r)}(t_i -x)
  K_r \lkr \frac{x-z_{l-1}}{\lamn} \rkr \right| dx \rkv^2 \nonumber\\
& \leq &  \sum_{i=1}^n   \| K_r \|^2_\infty \lkv  \sum_{l=1}^{M_n} \int_{ z_{l-1}}^{ z_{l}} | g^{(r)}(t_i -x) |  dx \rkv^2
\leq  n \,  \| g^{(r)} \|^2 \| K_r \|^2_\infty .
\label{Delta2}
\end{eqnarray}
Combining formulae \fr{eq:K_to_Q}--\fr{Delta2},  we obtain that, in
order to satisfy the condition  \fr{eq:cond2}, we need the following inequality to hold
\be \label{eq:inequa}
\KK (\PP_{f_1}, \PP_{f_2})  \leq
\frac{2 L^2 \lamn^{2m + 2r} n} {\sigma^2 T_n} \, \| K_r \|^2_\infty  [B_r^2 + \|g^{(r)} \|^2_2 ] \leq
\frac{1}{16}\, \frac{M_n \log 2}{8}.
\ee
Note that $\frac{T_n}{M_n} \lkr 1 - \frac{T_n}{M_n} \rkr \leq \lambda_n \leq \frac{T_n}{M_n}$.
Choosing  $M_n = C n^{1/(2(m+r)+1)} T_n^{(2(m+r)-1)/(2(m+r)+1)}$ and observing that
$T_n/M_n \rightarrow 0$ as $n \rightarrow \infty$,
obtain $\lamn \geq   T_n/(2 M_n)  \geq C  (T^2_n n^{-1})^{1/(2(m+r)+1)}$. Therefore,
both   conditions    \fr{eq:cond1} and \fr{eq:cond2} hold and theorem is proved.
\newline $\Box$
\\


\subsection*{Proof of Theorem \ref{th:semiexplicit}}
To prove Theorem \ref{th:semiexplicit} we use the following Lemma
\ref{lem:phi_properties} which can be viewed as  a version of  Theorem 7.2.4 of
Gripenberg, Londen and Staffans (1990, Chapter 7) adapted to our notations.

\begin{lemma}  \label{lem:phi_properties}
Let $s_g$ be such that
\be \label{cond_on_G}
\inf_{Re (s) = s_g} |\tilde{g}(s)| >0\ \ \ \mbox{and}\  \ \ \
\lim_{\stackrel{|s| \rightarrow \infty}{Re (s) \geq   s_g}} |s^r \tilde{g}(s)| > 0.
\ee
Then, solution $\phi(\cdot)$ of equation \fr{eq:new_eq} can be presented as
\be \label{phi_form2}
\phi(t) = \sum_{l=0}^L  \sum_{j=0}^{\alpha_l-1} \frac{a_{l,j}}{j!} t^j e^{s_l t} + \phi_1(t)
\ee
where $L$ is the total number of distinct zeros $s_l$ of $s^r \tilde{g}(s)$ such that $Re (s_l) > Re (s_g)$,
$\alpha_l$ is the order of zero $s_l$ and $\phi_1 \in L_1$.
\end{lemma}

\noindent
Choose $s_g$ such that $s^* < s_g <0$.
Then, the  first condition in \fr{cond_on_G} immediately follows from
Assumption (A2).
To validate the second assumption in \fr{cond_on_G}, note that
for $s = s_1 + i s_2$ conditions  $Re (s) \geq   s_g$ and  $|s| \rightarrow \infty$  imply that
either $s_1 \rightarrow \infty$ or $|s_2| \rightarrow \infty$, or both.
Recall  that $s^r \tilde{g}(s) =  B_r + \widetilde{G^{(r)}}(s)$.
If  $s_1 \rightarrow \infty$, no matter whether $s_2$ is finite or $s_2 \rightarrow \infty$,  one has
\be \label{s1_infinite}
\lim_{Re (s)   \rightarrow \infty} |s^r \tilde{g}(s)| =
\lim_{Re (s)   \rightarrow \infty} |B_r + \int_0^\infty g^{(r)}(t)e^{-st}dt| = |B_r| > 0.
\ee
If $s_1$ is finite, $s_1 \geq s_g$, and $|s_2| \rightarrow \infty$,  then Laplace transform
$\widetilde{g^{(r)}}(s) = \int_0^\infty g^{(r)} (t) e^{-st} dt$ is equal to Fourier transform ${\cal F} [g^{(r)} (t) e^{-s_1 t}] (s_2)$
of function $g^{(r)} (t) e^{-s_1 t}$ at the point $s_2$. Since  $g^{(r)} (t) e^{-s_1 t} \in L_1(\mathbb{R}^+)$, one obtains
$$
\lim_{|s_2|\rightarrow \infty} \int_0^\infty g^{(r)} (t) e^{-st} dt =
\lim_{|s_2|\rightarrow \infty}  {\cal F} [g^{(r)} (t) e^{-s_1 t}] (s_2) =0,
$$
and \fr{s1_infinite} holds again.
Hence, the second assumption in \fr{cond_on_G} is valid,
and Lemma \ref{lem:phi_properties} can be applied.

Note that, under Assumption (A2), $\tilde{g}(s)$ has no zeros with
$Re (s) > s_g$ and, therefore,  $s^r \tilde{g}(s)$ has a single zero of $r$-th order at $s=0$.
Lemma \ref{lem:phi_properties} yields then that
$\phi(t) =  \phi_0(t) +  \phi_1(t)$, where
\be \label{phio}
\phi_0(t)=\sum_{j=0}^{r-1} \frac{a_{0,j}}{j!} t^j,\;\;\;a_{0,j}=\phi^{(j)}(0),
\ee
and integrating by parts, one has
\be \label{int_by_parts}
\int_0^t q^{(r)}(t-\tau) \phi_0 (\tau)d\tau =
 \sum_{j=0}^{r-1} \phi_0^{(r-j-1)} (0) q^{(j)}(t),
\ee
that completes the proof of \fr{eq:f_semiexplicit}.

In order to prove \fr{eq:f_explicit} -- \fr{eq:bj}, note that it follows
from  equation \fr{eq:Lapl_phi} that $\tilde{\phi}(s)$ has poles $s_l$, $l=0, \ldots, M$, of
respective orders $\alpha_l$, where
$s_0=0$ and $\alpha_0 =r$. Since, by \fr{s1_infinite}, one has
$$
\lim_{\stackrel{|s| \rightarrow \infty}{Re (s) \geq   s_g}}
|s^r \tilde{g}(s)| > 0
$$
and, therefore, $\tilde{\phi}$ does not have a pole at infinity.
Then, $\tilde{\phi}$ is a rational function and, consequently,
can be represented using Cauchy integral formula
$$
\tilde{\phi}(s) = - \frac{1}{2\pi i} \sum_{l=0}^M \oint_{C_l}\frac{\tilde{\phi}(z)}{z-s} dz
$$
where $C_l$, $l=0,..., M$, is a circle around the pole $s_l$ such that this
circle does not enclose any other pole of $\tilde{\phi}$ (see LePage, 1961,
Section 5.14).
Using Laurent expansion of $\tilde{\phi}(z)$ around $s_l$, we have
\begin{eqnarray*}
I_l (s)  =    \frac{1}{2\pi i}  \oint_{C_l}\frac{\tilde{\phi}(z)}{z-s} dz
 = - \sum_{j=0}^{\alpha_l-1} \frac{1}{(s-s_l)^{j+1}}\ \frac{1}{(\alpha_l -1-j)!}\
\frac{d^{\alpha_l -j-1}}{ds^{\alpha_l -j-1}} \lkv (s - s_l)^{\alpha_l}
\tilde{\phi}(s) \rkv \Bigg|_{s=s_l}
\end{eqnarray*}
Combining the last two expressions and taking inverse Laplace transform of
$\tilde{\phi}(s)$ yields
$$
\phi(t) = \sum_{l=0}^M \sum_{j=0}^{\alpha_l-1} \frac{a_{l,j}}{j!} t^j e^{s_l t} = \phi_0 (t) + \phi_1(t),
$$
where $\phi_0$ is given by \fr{phio}, same as before, and
\begin{equation}  \label{eq:phi_one}
\phi_1 (t) =   \sum_{l=1}^M \sum_{j=0}^{\alpha_l-1} \frac{a_{l,j}}{j!} t^j e^{s_l t}.
\end{equation}
Repeat  calculations in \fr{int_by_parts} and also note that, by similar considerations,
for every $j = 0, ..., \alpha_l-1$, one can write
$$
\int_0^t q^{(r)}(t-x) x^j e^{s_l x} dx = \sum_{k=0}^{r-1} q^{(r-k-1)} (t)\,
\frac{d^k}{dx^k} \lkv x^j e^{s_l x} \rkv \Bigg|_{x=0} +
\int_0^t q(t-x) \frac{d^{(r-1)}}{dx^{(r-1)}} \lkv x^j e^{s_l x} \rkv dx.
$$
To complete the proof,  evaluate the derivatives, observe that
$$
\frac{d^k}{dx^k} \lkv x^j e^{s_l x} \rkv \Bigg|_{x=0} = {k \choose j}\, s_l^{k-j}
$$
and interchange summation with respect to $j$ and $k$.
\newline $\Box$


\subsection*{Proof of Theorem \ref{th:risk_upper_bounds}}

Since the estimator  \fr{eq:f_explicit_est} is just a particular form of the
estimator \fr{eq:estimator}, it is sufficient to carry out the proof for the estimator
\fr{eq:estimator} of $f$. From \fr{eq:f_semiexplicit},  one immediately obtains
\begin{eqnarray}
E||\hat{f}_n-f||^2_{[0,T_n]} \leq \frac{r+2}{B^2_r}
\left( E||\widehat{q_{\hat{\lambda}_{r,n}}^{(r)}}-q^{(r)}||^2_{[0,T_n]} \right. & + &
\sum_{j=0}^{r-1}a^2_{0,r-1-j}E||\widehat{q_{\hat{\lambda}_{j,n}}^{(j)}}-q^{(j)}||^2_{[0,T_n]}  \nonumber\\
& + &
\left. ||\widehat{q_{\hat{\lambda}_{r,n}}^{(r)}} \ast \phi_1 -  q^{(r)} \ast \phi_1||^2_{[0,T_n]}\right), \label{eq:mse}
\end{eqnarray}
where $\widehat{q_{\hat{\lambda}_{j,n}}^{(j)}},\;j=0,\ldots,r$ are given in \fr{eq:deriv_est}.

The proof is based on the  following proposition which provides  upper bounds for the risks
$E||\widehat{q_{\hat{\lambda}_{j,n}}^{(j)}}-q^{(j)}||^2_{[0,T_n]},\;j=0,...,r$ in \fr{eq:mse}.

\begin{proposition} \label{prop:deriv_risk}
Let condition \fr{k_cond} and  Assumptions (A1)-(A4) hold. Let kernel $K_j$
be of order $(L,j)$, where $L > r$ and $0 \leq j \leq r$, and satisfies Assumptions
(K1) and (K2).  Then, for all $A'>0$,
\be \label{eq:LepDeriv}
\sup_{q  \in W^{m+r}(A')}\, E|| \widehat{q_{\hat{\lambda}_{j,n}}^{(j)}} - q^{(j)}||^2_{[0,T_n]}  =
O \lkr \lkr \frac{T_n^2}{n} \rkr^\frac{2(r+m-j)}{2(r+m)+1} \rkr.
\ee
\end{proposition}
In particular, Proposition \ref{prop:deriv_risk} implies that the errors of
estimating $q^{(j)}$ in \fr{eq:mse} are dominated by the estimation error
of the highest order derivative $q^{(r)}$. Furthermore,
$\phi_1 \in L_1$ (see Theorem \ref{th:semiexplicit}) and, therefore,
$$
\|\widehat{q_{\hat{\lambda}_{r,n}}^{(r)}}\ast \phi_1 -  q^{(r)} \ast \phi_1  \|_2 \leq \|\phi_1\|_1 \cdot
\|\widehat{q_n^{(r)}}-q^{(r)}\|_2 = O\left(\|\widehat{q_{\hat{\lambda}_{r,n}}^{(r)}}-q^{(r)}\|_2\right)
$$
(see also Theorem 2.2.2 of Gripenberg, Londen and Staffans, 1990).

Thus, \fr{eq:mse} and Proposition \ref{prop:deriv_risk} yield
$$
E||\hat{f}_n-f||^2_{[0,T_n]}=O\left(E||\widehat{q_{\hat{\lambda}_{r,n}}^{(r)}}-q^{(r)}||^2_{[0,T_n]}
\right)=O\lkr \lkr \frac{T_n^2}{n} \rkr ^\frac{2m}{2m+2r+1} \rkr
$$
\newline $\Box$


\subsection*{Proof of Proposition \ref{prop:deriv_risk}}

For simplicity of notations we drop the index $n$ in $\hat{\lambda}_{j,n}$.

Recall that under Assumptions (A1)-(A3), $q \in W^{r+m}$ (see Section
\ref{subsec:explicit}).
By the standard asymptotic calculus for kernel estimation (see, e.g.,
Gasser and M\"uller, 1984) for estimator \fr{eq:pc} and any  interior point $t$ of $(0, T_n)$, one then has
$$
Var\left(\widehat{q_{\lambda}^{(j)}}(t)\right)=
\frac{\sigma^2}{\lambda_j^{2(j+1)}} \ \sum_{i=1}^n (t_i-t_{i-1})^2 K_j^2 \lkr \frac{t_i-t}{\lambda_j} \rkr =
\frac{\sigma^2}{\lambda_j^{2j+1}} \frac{T_n}{n}\int K_j^{2}(u)du \; (1+o(1)).
$$
The required boundary corrections ensure
the same order of error for the values  $t$ close to the boundaries
(Gasser and M\"uller, 1984) and the integrated variance then is
\be \label{eq:var_comp}
V_j(\lambda_j) = \int_0^{T_n} Var\left(\widehat{q_{\lambda_j}^{(j)}}(t)\right)dt =
V_{0j}\ \frac{T^2_n}{\lambda_j^{2j+1}n}\;(1+o(1)),
\ee
where $V_{0j} = \sigma^2 ||K_j||^2$.
Similarly,   the integrated squared bias can be written   as
\be \label{eq:bias_comp}
B_j^2(\lambda_j, q)=\int_0^{T_n}\left(E\left(\widehat{q_{\lambda_j}^{(j)}}(t)\right)-q^{(j)}(t)\right)^2dt =
B_{0j}  \lambda_j^{2(r+m-j)}(1+o(1)),
\ee
where $B_{0j} =B^{-1}_0 ||q^{(r+m)}||^2 ||K_j||^2$ and $B_0=2\left((r+m-1)!\right)^2
(2(r+m)-1)(2(r+m)+1)$.
Hence,
\be \label{eq:band_calc}
\sup_{q \in W^{m+r}(A)} E||\widehat{q_{\lambda_j}^{(j)}}-q^{(j)}||^2_{L_2([0,T_n)]}=
\sup_{q \in W^{m+r}(A)} \left(V_j(\lambda_j)+B_j^2(\lambda_j, q)\right)
=O\left(\frac{T^2_n}{\lambda_j^{2j+1}n}\right)+O\left(\lambda_j^{2(r+m-j)}\right).
\ee
It follows from \fr{eq:var_comp}  and \fr{eq:bias_comp} that the
asymptotically optimal bandwidth that minimizes
$E||\widehat{q_{\lambda_j}^{(j)}}-q^{(j)}||^2_{L_2([0,T_n])}$ is
\be \label{eq:optband}
\lambda_j^*=O\left(\left(\frac{T_n^2}{n}\right)^{\frac{1}{2(r+m)+1}}\right)
\ee
and the corresponding risk of estimating $q^{(j)}$ is given by
\be \label{eq:mseglobal}
\sup_{q \in W^{m+r}(A')} E||\widehat{q_{\lambda_j^*}^{(j)}}-q^{(j)}||^2_{L_2([0,T_n])}=
O\left(\left(\frac{T_n^2}{n}\right)^{\frac{2(r+m-j)}{2(r+m)+1}}\right).
\ee

Now we need to prove that \fr{eq:mseglobal} remains valid
when $\lambda_j^*$  is replaced by $\hat{\lambda}_j$ selected by Lepski procedure,
that is,
$$
\sup_{q \in W^{m+r}(A')}\, E||\widehat{q_{\hat{\lambda}_j}^{(j)}} - q^{(j)}||^2  =
O \lkr \lkr \frac{T_n^2}{n} \rkr^\frac{2(r+m-j)}{2(r+m)+1} \rkr
$$
for all $A'>0$.
Set $d_j$ and  $\lambda^*_j$ in \fr{eq:optband} to be, respectively,
$$
d_j = \frac{C_j- \mu ||K_j||}{2 ||K_j||}, \quad
\lambda^*_j = \lkr d_j^2\ \frac{\sigma^2 B_0}{2(A')^2 }\,
\frac{T_n^2}{n} \rkr^{\frac{1}{2(r+m)+1}},
$$
where $C_j$ is defined in \fr{c0}.
Note that
$$
E \| \widehat{q_{ \hat{\lambda}_j}^{(j)}} - q^{(j)} \|^2  =
E \left\{ \| \widehat{q_{ \hat{\lambda}_j}^{(j)}} - q^{(j)} \|^2 I(\hat{\lambda}_j \geq \lambda^*_j) \right\}
+
E \left\{ \| \widehat{q_{ \hat{\lambda}_j}^{(j)}} - q^{(j)} \|^2 I(\hat{\lambda}_j < \lambda^*_j) \right\}
= \Delta_1 + \Delta_2.
$$
For $\hat{\lambda}_j \geq \lambda^*_j$, equations \fr{eq:mseglobal} and \fr{eq:adapband}
imply that uniformly over $q \in W^{m+r}(A')$
\begin{eqnarray}
\Delta_1 & \leq & 2 E \left\{ \| \widehat{q_{ \hat{\lambda}_j}^{(j)}} - q_{{\lambda}^*_j}^{(j)} \|^2 I(\hat{\lambda}_j > \lambda^*_j) \right\}
+ 2 E \left\{ \| \widehat{q_{ \hat{\lambda}^*_j}^{(j)}} - q^{(j)} \|^2 I(\hat{\lambda}_j > \lambda^*_j) \right\}\nonumber\\
& = &  O \lkr n^{-1} T_n^2 (\lambda^*_j)^{-(2j+1)} \rkr +
O\lkr(n^{-1} T_n^2)^{-\frac{2(r+m-j)}{2(r+m)+1}} \rkr
=  O \lkr (n^{-1} T_n^2)^{-\frac{2(r+m-j)}{2(r+m)+1}} \rkr.  \label{delta1}
\end{eqnarray}
For $(n^{-1} T_n^2)^{\frac{1}{2j+1}} \leq \hat{\lambda}_j < \lambda^*_j$,
by direct calculus similar to that carried out above, one can show that
\begin{eqnarray*}
\sup_{q   \in W^{m+r}(A)} \, E   \| \widehat{q_{ \hat{\lambda}_j}^{(j)}} - q^{(j)} \|^4   & = &
O \lkr  (\lambda^*_j)^{-2(2j+1)} n^{-2} T_n^4 \rkr + O((\lambda^*_n)^{4(r+m-j)}) = O(1).
\end{eqnarray*}
Hence,
\begin{eqnarray}
\sup_{q   \in W^{m+r}(A)} \Delta_2 & \leq &
\sup_{q   \in W^{m+r}(A)} \sqrt{E   \| \widehat{q_{ \hat{\lambda}_j}^{(j)}} - q^{(j)} \|^4 }\ \sqrt{ P(\hat{\lambda}_j < \lambda^*_j)}
\nonumber\\
& = & O \lkr \sup_{q \in W^{m+r}(A)} \sqrt{ P(\hat{\lambda}_j < \lambda^*_j)}  \rkr.  \label{delta2}
\end{eqnarray}

If $\lambda^*_j > \hat{\lambda}_j$, it follows from  definition \fr{eq:adapband} of $\hat{\lambda}_j$
that there exists $\tilde{h}  < \lambda^*_j$ such that
$
\| \widehat{q_{\lambda^*_j}^{(j)}} - \widehat{q_{\tilde{h}}^{(j)}} \|^2 >  4\, C_j^2 n^{-1}\sigma^2 T_n^2 \tilde{h}^{-(2j+1)},
$
where, by  \fr{c0} and definition of $d_j$, we have  $C_j = \|K_j\| (\mu + 2 d_j)$.
It follows from \fr{eq:var_comp} and \fr{eq:bias_comp} that, for all
$h < \lambda^*_j$, the variance term dominates  the squared bias, that is,
$$
\sup_{q   \in W^{m+r}(A')}\, \|E\widehat{q_{h}^{(j)}}-q^{(j)}\|^2 \leq
d_j^2 \sigma^2 \|K_j\|^2 n^{-1}T_n^2 h^{-(2j+1)}.
$$
Hence, for all $\tilde{h} < \lambda^*_j$ and $q \in W^{m+r}(A')$, one has
\begin{eqnarray*} 
 P \lkr \| \widehat{q_{\lambda^*_j}^{(j)}} - \widehat{q_{\tilde{h}}^{(j)}} \|^2 >  4\, C_j^2 n^{-1}\sigma^2 T_n^2 \tilde{h}^{-(2j+1)} \rkr
& < &
P \lkr \| \widehat{q_{\lambda^*_j}^{(j)}}  - E q_{\lambda^*_j}^{(j)} \|^2 > \sigma^2 \|K_j\|^2 (\mu +   d_j)^2 n^{-1} T_n^2 \tilde{h}^{-(2j+1)} \rkr \\
& + &
P \lkr \| \widehat{q_{\tilde{h}}^{(j)}}  - E q_{\tilde{h}}^{(j)} \|^2 > \sigma^2 \|K_j\|^2 (\mu +   d_j)^2 n^{-1} T_n^2 \tilde{h}^{-(2j+1)} \rkr \\
\end{eqnarray*}
due to $C_j  - \|K_j\| d_j  > \|K_j\|  (\mu + d_j)$.
Thus, uniformly over $q \in W^{s+r}(A')$, one has
\begin{eqnarray}
P(\hat{\lambda}_j < \lambda^*_j) & \leq & \sum_{\stackrel{h \in \Lambda_j}{h \leq \lambda^*_j}}
P(\tilde{h}=h)\ P \lkr \| \widehat{q_{\lambda^*_j}^{(j)}}  - \widehat{q_h^{(j)}} \|^2 >
4\, \sigma^2 C_j^2  n^{-1} T_n^2 h^{-(2j+1)} \rkr \label{eq:P_sum}
\\
& \leq & 2\ \sum_{\stackrel{h \in \Lambda_j}{h \leq \lambda^*_j}} P(\tilde{h}=h)\
P \lkr \| \widehat{q_h^{(j)}} - E \widehat{q_h^{(j)}} \|^2 \geq
\sigma^2   \|K_j\|^2 (\mu + d_j)^2   n^{-1} T^2_n h^{-(2j+1)} \rkr. \nonumber
\end{eqnarray}
Note that
$$
 \| \widehat{q_h^{(j)}} - E \widehat{q_h^{(j)}} \|^2  =
\left\|\sum_{i=1}^n h^{-(j+1)} K_j \lkr \frac{t-t_i}{h} \rkr (t_i-t_{i-1}) \epsilon_i \right\|^2
= h^{-(2j+1)} n^{-2} T_n^2\ \bepsilon^T \bQ \bepsilon,
$$
where $\bQ$ is an $n \times n$ symmetric nonnegative-definite matrix with elements
\be \label{eq:Q}
Q_{il} =  \frac{n^2}{T_n^2} (t_i - t_{i-1}) (t_l - t_{l-1})\, \int_{-1}^1 K_j(z)  K_j \lkr z + \frac{t_i - t_l}{h} \rkr  dz.
\ee
Then,
\begin{equation} \label{eq:chisq}
P \lkr \| \widehat{q_h^{(j)}} - E \widehat{q_h^{(j)}} \|^2 \geq
\sigma^2   \|K_j\|^2 (\mu + d_j)^2   n^{-1} T^2_n h^{-(2j+1)}  \rkr =
P\lkr \bepsilon^T \bQ \bepsilon \geq   n\, \sigma^2   \|K_j\|^2 (\mu + d_j)^2 \rkr.
\end{equation}
Applying a $\chi^2$-type inequality which initially appeared  in  Laurent
and Massart (1998), was improved by Comte (2001) and furthermore
by Gendre (2013), we derive that, for any $x>0$,
\begin{equation} \label{eq:chisq1}
P \lkr \sigma^{-2} \bepsilon^T \bQ \bepsilon \geq \lkv \sqrt{\Tr(\bQ)} +
 \sqrt{ x \rho_{\max}^2(\bQ)} \rkv^2   \rkr\leq e^{-x},
\end{equation}
where $\Tr(\bQ)$ is the trace of $\bQ$,  and $\rho_{\max}^2(\bQ)$ is the maximal eigenvalue of $\bQ$.
Note that
$$
\Tr (\bQ) = \frac{n^2}{T_n^2} \sum_{i=1}^n (t_i - t_{i-1})^2 \|K_j \|^2 \leq n \mu^2 \|K_j\|^2.
$$
and $\rho_{\max}^2(\bQ)$ is the spectral norm of matrix $\bQ$  which is dominated by
any other norm. In particular,
$$
\rho_{\max}^2(\bQ) \leq  \max_k \sum_{l=1}^n |Q_{kl}|
=\frac{n^2}{T_n^2} \max_k (t_k - t_{k-1}) \int_{-1}^1 |K_j(z)| \lkv  \sum_{l=1}^n
\left| K_j\lkr z + \frac{t_k - t_l}{h} \rkr \right| (t_l - t_{l-1}) \rkv dz.
$$
Since
\begin{eqnarray*}
\sum_{l=1}^n \left|K_j \lkr z + \frac{t_k - t_l}{h} \rkr \right| (t_l - t_{l-1})
& = & \int_{-1}^1 \left|K_j \lkr z + \frac{t_k - t}{h} \rkr \right| dt (1+o(1)) \\
& = & h \int_{-1}^1  \left|K_j \lkr z +  \frac{t_k}{h} - y \rkr \right| dt (1+o(1)),
\end{eqnarray*}
we derive
\begin{eqnarray*}
\rho_{\max}^2(\bQ)  & \leq  &   \frac{n^2}{T_n^2} \max_k \lkv (t_k - t_{k-1})\, h \,
\int_{-1}^1 \int_{-1}^1 |K_j(z)| |K_j(z +  t_k/h - y)| dz dy \rkv \\
& \leq &   \mu \frac{nh}{T_n} \lkv \int_{-1}^1 |K_j(z)|dz \rkv^2 \leq
 2 \mu \|K_j\|^2\, \frac{nh}{T_n}.
\end{eqnarray*}

Using inequality \fr{eq:chisq1} with $x= d_j^2 T_n/ (2\mu h)$ and $h < \lambda^*_j$
one obtains
\be \label{eq:expon_ineq}
P\lkr \| \widehat{q_h^{(j)}} - E \widehat{q_h}^{(j)} \|^2
\geq \frac{\sigma^2   \|K_j\|^2 (\mu + d_j)^2   T^2_n}{n h^{2j+1}} \rkr
\leq \exp \lkr - \frac{d_j^2 T_n}{2\mu h} \rkr \leq
\exp \lkr - c_j n^{\frac{1}{2(r+m)+1}}
T_n^{\frac{2(r+m)-1}{2(r+m)+1}} \rkr
\ee
where $c_j$ depends on $m$, $A'$, $\mu$ and $d_j$.
Combination of \fr{delta1}, \fr{delta2}, \fr{eq:P_sum} and \fr{eq:expon_ineq} completes the proof.
\newline $\Box$
\\


{\bf Proof of Lemma  \ref{th:suppl_lemma}.   }
Definitions \fr{K_definition} imply that $k (x) = K_1^\prime (x)$,
$K_{j-1}^\prime (x) = K_j (x)$ and $K_j (0) =0$, $j=1,...,r$.
Observe that condition $K_j (1) =0$, $j=1,...,r$, is equivalent to
\be \label{suppl}
\int_0^1 K_j(x) dx = 0, \ \ \ j=0, ..., r-1,
\ee
where $K_0(x) = k(x)$.  It is easy to see that   \fr{suppl} is valid for $j=0$.
For $j \geq     1$, note that, by formula (4.631) of Gradshtein and  Ryzhik (1980),
\be \label{expr}
K_j(x) = \int_0^x dz_{j-1} \int_0^{z_{j-1}} dz_{j-2}  \ldots \int_0^{z_1} k(z) dz =
\frac{1}{(j-1)!} \ \int_0^x (x-z)^{j-1} k(z) dz.
\ee
Then, for any $x \in [0,1]$, one has
$|K_j(x)| \leq  [(j-1)!]^{-1}\, \|k\|_\infty\,  \int_0^x (x-z)^{j-1} dz \leq  \|k\|_\infty$.
Moreover, by \fr{expr}, for $j=1, \ldots, r-1$,  one has
\begin{align*}
& \int_0^1 K_j(x) dx = \frac{1}{(j-1)!} \, \int_0^1 dx \int_0^x (x-z)^{j-1} k(z) dz \\
& = \frac{1}{(j-1)!} \, \int_0^1 k(z) dz \int_z^1 (x-z)^{j-1} dx =
\frac{1}{(j-1)! j!} \, \int_0^1 (1-z)^j\, k(z) dz =0.
\end{align*}

Now, it remains to prove formula \fr{int_bound}. Note that support of the function
$k(u/\lamn - (l-1))$ coincides with $(z_{l-1}, z_l)$,  so that
\be \label{form1}
I(i,l) = \int_0^{t_i} g(t_i -x) k \lkr \frac{x-z_{l-1}}{\lamn} \rkr dx =
\int_{\min(z_{l-1}, t_i)}^{\min(z_{l}, t_i)} g(t_i -x) k \lkr \frac{x-z_{l-1}}{\lamn} \rkr dx.
\ee
Formula \fr{form1}  implies that $I(i,l) =0$ whenever $z_{l-1} \geq y_i$. If $z_{l-1} < y_i \leq z_l$,
it follows from \fr{form1} that
$$
I(i,l) = \int_{z_{l-1}}^{t_i} g(t_i -x) k \lkr \frac{x-z_{l-1}}{\lamn} \rkr dx.
$$
Introduce new variable $t= x-z_{l-1}$ and denote  $u_{il} = t_i - z_{l-1}$. Then, recalling condition \fr{k_cond}
and using integration by parts, we derive
\begin{eqnarray*}
I(i,l) & = &  \int_0^{u_{il}} g(u_{il}-t) k \lkr \frac{t}{\lamn} \rkr dt =
\lamn g(u_{il}-t) K_1 \lkr \frac{t}{\lamn} \rkr \Bigg|_0^{u_{il}} + \lamn \int_0^{u_{il}} g^\prime (u_{il}-t)  K_1 \lkr \frac{t}{\lamn} \rkr dt \\
& = &  \ldots =  \lamn^r g^{(r-1)} (u_{il}-t) K_r \lkr \frac{t}{\lamn} \rkr \Bigg|_0^{u_{il}}
+ \lamn^r \int_0^{u_{il}} g^{r} (u_{il}-t)  K_r \lkr \frac{t}{\lamn} \rkr dt.\\
\end{eqnarray*}
Changing variables back to $x$, we arrive at
\be \label{form2}
I(i,l) = \lamn^r \lkv B_r\, K_r \lkr \frac{t_i - z_{l-1}}{\lamn} \rkr +
\int_{z_{l-1}}^{t_i}   g^{(r)}(t_i -x) K_r \lkr \frac{x-z_{l-1}}{\lamn} \rkr dx \rkv.
\ee

Finally, consider the case when $z_l \leq y_i$. Then, using relation  $z_l= z_{l-1} + \lamn$,
integration by parts and the fact that $K_j (0) = K_j(1) =0$ for $j=1,...,r$, we obtain
\begin{eqnarray*}
I(i,l) & = &   \int_{z_{l-1}}^{z_l} g(t_i -x) k \lkr \frac{x-z_{l-1}}{\lamn} \rkr dx
= \lamn \int_0^1 g(t_i - z_{l-1} - \lamn t ) k(t) dt \\
& = &   \ldots =  \lamn^{r+1}  \int_0^1 g^{r} (t_i - z_{l-1} - \lamn t )  K_r (t) dt
= \lamn^r  \int_{z_{l-1}}^{z_l} g^{(r)} (t_i -x) K_r \lkr \frac{x-z_{l-1}}{\lamn} \rkr dx
\end{eqnarray*}
which, in combination with \fr{form2}, completes the proof.
\newline $\Box$
\\

\noindent
Felix Abramovich\\
Department of Statistics $\&$ Operations Research \\
Tel Aviv University \\
Tel Aviv 69978,  Israel \\
{\em felix@post.tau.ac.il} \\

\noindent
Marianna Pensky\\
Department of Mathematics \\
University of Central Florida \\
Orlando FL 32816-1353, USA \\
{\em Marianna.Pensky@ucf.edu}\\

\noindent
Yves Rozenholc\\
Universit\'{e} Paris Descartes\\
MAP5-UMR CNRS 8145\\
75270 Paris Cedex, France\\ 
{\em yves.rozenholc@univ-paris5.fr}

\end{document}